\documentstyle[12pt,leqno,amsfonts,amscd]{amsart}
\newtheorem{Theorem}{Theorem}[section]

\newtheorem{Proposition}[Theorem]{Proposition}

\newtheorem{Lemma}[Theorem]{Lemma}
\newtheorem{Corollary}[Theorem]{Corollary}
\theoremstyle{remark}

\begin{document}
\title{Transcendence measures and algebraic growth of entire
functions}
\date{March 2004}
\author{ Dan Coman and Evgeny A. Poletsky}
\thanks{Both authors are supported by NSF Grants.}
\address{ Department of Mathematics,  215 Carnegie Hall,
Syracuse University,  Syracuse, NY 13244-1150, USA. E-mail:
dcoman@@syr.edu, eapolets@@syr.edu}
\begin{abstract}In this paper we obtain estimates for certain
transcendence measures of an entire function $f$. Using these
estimates, we prove Bernstein, doubling and Markov inequalities
for a polynomial $P(z,w)$ in ${\Bbb C}^2$ along the graph of $f$.
These inequalities provide, in turn, estimates for the number of
zeros of the function $P(z,f(z))$ in the disk of radius $r$, in
terms of the degree of $P$ and of $r$.
\par Our estimates hold for arbitrary entire functions $f$ of finite
order, and for a subsequence $\{n_j\}$ of degrees of polynomials.
But for special classes of functions, including the Riemann
$\zeta$-function, they hold for all degrees and are asymptotically
best possible. From this theory we derive lower estimates for a
certain algebraic measure of a set of values $f(E)$, in terms of
the size of the set $E$.
\end{abstract}
\maketitle
\section{Introduction}
\par In recent years there was a significant interest in the
behavior of a polynomial $P$ along an algebraic subvariety $X$ of
${\Bbb R}^n$ or ${\Bbb C}^n$. This started with the paper
\cite{FN} of Fefferman and Narasimhan, where they obtained local
{\it doubling inequalities}, which bound the ratio of the uniform
norms of $P$ on two concentric balls in $X$, in terms of the
degrees of $P$ and $X$, and of the ratio of the radii of these
balls.
\par Later, these inequalities were improved in papers of Brudnyi
\cite{Br} and Roytwarf and Yomdin \cite{RY},  and they were
applied to questions from analytic geometry, pseudodifferential
operators, to Hilbert's 16th problem, and so on.
\par Much earlier, Tijdeman \cite{Ti1} studied the behavior of
a polynomial $P(z,w)$ in ${\Bbb C}^2$ along the graph of the
exponential function $w=e^z$. In this situation, he obtained
global doubling inequalities and estimates for the number of zeros
of the function $P(z,e^z)$ in a disk of radius $r$. He used these
results in \cite{Ti2} to get new advancements in transcendental
number theory. The proofs in \cite{Ti2} involved transcendence
measures of numbers, which were studied extensively in
transcendental number theory.
\par Transcendence measures appear quite naturally when transcendental
objects are investigated. In general, if $B$ is a subring of a
commutative ring $A$, then an element $\omega\in A$ is called {\it
transcendental} over $B$ if $P(\omega)\ne0$, for any non-trivial
polynomial $P\in B[x]$. For example, if $A={\Bbb C}$ and $B={\Bbb
Z}$ we get the transcendental numbers, and if $A$ is the ring of
entire functions and ${\Bbb C}[z]$ the ring of polynomials in
${\Bbb C}$, we get the entire transcendental functions.
\par If $A$ and $B$ are normed rings and the algebra $B[x]$ is graded,
i.e., there is an increasing sequence of sets $B_n[x]$ such that
$\bigcup_{n\ge0}B_n[x]=B[x]$, then we can measure the
transcendence of $\omega$. For this, we define a suitable norm
$h(P)$ of $P\in B[x]$, and let  the {\it transcendence measure}
$\tau(\omega,n,H)$ of $\omega$ be the infimum of $\|P(\omega)\|$
over all polynomials $P\in B_n[x]$ with $1\le h(P)\le H$.
\par In our papers \cite{CP1} and \cite{CP2}, we started to study
{\it transcendence measures} of an entire function $f$. A
transcendence measure can be defined by
$$E_n(f)=\sup\{\|P\|_{\Delta^2}:\,P\in{\Bbb C}[z,w],\,\deg P\le
n,\,\|P(z,f(z))\|_\Delta\le1\}.$$ Here $\Delta$ is the closed unit
disk in ${\Bbb C}$ and the norms are uniform norms. Since $f$ is
usually fixed, we write $E_n=E_n(f)$, and let $e_n=\log E_n$. In
\cite{CP1} we proved that
$$e_n(e^z)=\frac12\,n^2\log n+O(n^2).$$
\par This transcendence measure is closely connected
with the following aspects of analysis and geometry:
\par 1) Polynomial estimates on ${\Bbb C}^2$: if $P(z,w)$ is a
polynomial of degree $n$ and $|P(z,f(z))|\le1$ on the unit disk
$\Delta$, then
$$|P(z,w)|\le E_n\exp\left(n\max\{\log^+|z|,\log^+|w|\}\right);$$
\par 2) Polynomial estimates along the graph of $f$: if
$$m(r)=m(r,f)=\max\{\log^+|f(z)|:\,|z|=r\},$$ $P_f(z)=P(z,f(z))$,
$\Delta_r=\{z\in{\Bbb C}:\,|z|\le r\}$ and
$$m_n(r)=m_n(r,f)=
\sup\{\log\|P_f\|_{\Delta_r}:\,\deg P\le n,\,
\|P_f\|_\Delta\le1\},$$ then
$$m_n(r)\le\frac{2e_n}{\log t_n}\,\log r,\;1\le r\le t_n,$$
where $t_n$ is defined by $nm(t_n)=e_n$. The functions $m_n(r,f)$,
$r>1$, can also be considered as transcendence measures of $f$, by
using $\|P_f\|_{\Delta_r}$ as the norm of a polynomial $P(z,w)$.
\par 3) Estimates on the number of zeros: if $Z_n(r)=Z_n(r,f)$ is the
maximum number of zeros of the function $P_f$ in the disk
$\Delta_r$ when $\deg P\le n$, then  $Z_n(r)\le2m_n(3r)$. The
number $Z_n(r)$ gives the maximum number of intersection points
of an algebraic variety of degree $n$ with the graph of $f$ in
${\Bbb C}^2$ lying over $\Delta_r$.
\par These connections were proved in \cite{CP2}, where we also
found an approach to estimate $e_n$ for general transcendental
functions. It allowed us to handle, in particular, the class of
functions $f(z)=e^{P(z)}$, where $P$ is a polynomial.
\par For any transcendental function $f$ one has
\cite[Proposition 1.3]{CP1}
$$m_n(r,f)\geq\frac{n^2+3n}{2}\,\log r,\;r\geq 1.$$
Using the transcendence measure $m_n(e,f)$, we define the {\it
lower order of transcendence} as
$$\underline\tau(f)=\sup\left\{\tau:\,
\liminf_{n\to\infty}\frac{m_n(e,f)}{n^\tau}>0\right\},$$ and the
{\it upper order of transcendence} as
$$\overline\tau(f)=\inf\left\{\tau:\,
\limsup_{n\to\infty}\frac{m_n(e,f)}{n^\tau}<\infty\right\}.$$
\par Since $m_n(e,f)\geq n^2/2$, we have $\underline\tau(f)\ge2$.
If $f(z)=e^z$ then $\lim_{n\rightarrow\infty}m_n(e)/n^2=1/2$
\cite[Theorem 1.2]{CP1}. More generally, if $f(z)=e^{P(z)}$, for
some polynomial $P$, then $\underline\tau(f)=\overline\tau(f)=2$
\cite[Theorem 5.1]{CP2}. For $\tau\geq3$, we constructed examples
of entire functions of order 1 and type $1/e$ with
$\tau-1\le\overline\tau(f)\le\tau$ \cite[Corollary 6.2]{CP2}. In
all these examples, $\underline\tau(f)=2$. Whether this was true
in general remained unsettled until the present paper.
\par The approach in \cite{CP2} was based on estimates of
$e_n(f)$ in terms of the $n$-th diameter of the set of preimages
of a point on the unit circle. The $n$-th diameter $d_n(F)$ of a
set $F$ is the minimal sum of radii of $n$ disks covering $F$. In
\cite{So}, using the theory of Dufresnoy, Sodin gave lower bounds
for the smallest number of disks of radius $R^\alpha$, $\alpha<1$,
needed to cover the set $f^{-1}(\{0,1\})\cap\Delta_R$, when $f$ is
a function of finite positive order $\rho$. Applied to our
problem, his result leads to only exponential estimates for $e_n$.
In Section \ref{S:end}, using the Ahlfors theory of covering
surfaces and certain results of Dufresnoy, we obtain the necessary
estimates for the $n$-th diameter. The results we need from these
theories are recalled in Section \ref{S:bf}.
\par The estimates for the $n$-th diameter allow us to obtain several
results, which can be summarized in the following theorem. In this
theorem, the second inequality is usually called a {\it Bernstein
} inequality, the third -- a {\it Bezout} inequality, the fourth
-- a {\it doubling} inequality, and the fifth -- a {\it Markov}
inequality. Bernstein and Markov inequalities have been
extensively studied and have wide applications, for example in
approximation theory (see e.g. \cite{BBLT} and references
therein).
\begin{Theorem}\label{T:ta} For any entire function $f$ of finite order
$\rho>0$, there exist sequences of integers $\{n_j\}$ and
$\epsilon_j>0,\;\epsilon_j\to0$, such that
$$e_{n_j}\leq
C_1n_j^2\log n_j\,,\;\;m_{n_j}(r)\leq C_2n_j^2\log r,\;1\leq
r\leq\frac{1}{2}\,n_j^{1/\rho-\epsilon_j}.$$ For every $r\geq1$
there exists an integer $j_r$ such that if $j\ge j_r$ then
$$Z_{n_j}(r)\le C_3n_j^2\,,\;\;\frac{M(2r,P_f)}{M(r,P_f)}\le 2^{an_j^2}\,,
\;\;M(r,P_f')\le C_4n_j^2\frac{M(r,P_f)}{r}\,,$$ where $P$ is a
polynomial of degree at most $n_j$.
\end{Theorem}
Here $M(r,F)=\max\{|F(z)|:\,|z|=r\}$ and the constants are
effectively computed and depend only on $\rho$. A sequence of
integers $\{n_j\}$ for which the above theorem holds will be
called a {\it fundamental sequence} for $f$. It follows from this
theorem that $\underline\tau(f)=2$ for all entire functions of
finite positive order.
\par Theorem \ref{T:ta} is proved in Sections \ref{S:gpe}, \ref{S:loot}
and \ref{S:di}, where we also show that for entire functions with
a covering system of admissible intervals
$I(R,\alpha,\beta,\gamma,C)$ (see Section  \ref{S:gpe}), the
inequalities in Theorem \ref{T:ta} hold for {\it all} $n$
sufficiently large. Again, all constants are effectively computed.
The only change is that one should substitute $n^{1+1/\gamma}$
instead of $n_j^2$.
\par In Section \ref{S:scof} we give three sufficient criteria for
classes of functions to have a covering system of admissible
intervals $I(R,\alpha,\beta,\gamma,C)$. The first one states that
if $A_1m(r,f)\le m(kr,f)\le A_2m(r,f)$, for some constants
$A_1,A_2,k>1$, then the function $f$ has a covering system of
admissible intervals $I(R,\alpha,\beta,1,C)$. This class includes
all functions $f(z)=\sum_{j=1}^mp_j(z)e^{q_j(z)}$, where $p_j$ and
$q_j$ are polynomials, and, as shown in Section \ref{S:rf}, the
Riemann $\zeta$-function and the function $\xi$. It follows that
for such functions Theorem \ref{T:ta} holds for all $n$
sufficiently large.
\par The second criterion can be applied when we know that
$m(r,f)\le r^{\phi(r)}$ and $r^{\phi(r)-\rho}$ is a slowly
increasing function (see Theorem \ref{T:pesif}). Finally,
Corollary \ref{C:tc} gives a criterion based on the behavior of
the Taylor coefficients of $f$, similar to the formulas for the
order and type of $f$.
\par In Section \ref{S:ext} we introduce and study an {\it extremal
function} $W^\star(z)$, related to Bernstein inequalities, and we
prove that $W^\star(z)=\frac12\,\log^+|z|$ when $f(z)=e^z$.
\par In Section \ref{S:agtf} we address a problem
posed by Mahler in \cite{M}: given an entire transcendental
function $f$, describe, or at least find properties of, the set of
algebraic numbers where the values of $f$ are also algebraic.
There are many results claiming that this set is finite when
either $f$ is a special function, or when {\it all} the
derivatives of $f$ take algebraic values on this set and their
algebraic measure satisfies some growth conditions (see, e.g.,
\cite{Sc}, \cite{St}, \cite{La}, \cite{W}). But a general entire
function may take algebraic values on any set of algebraic numbers
(see \cite{M} and \cite{GS}), in particular, on any algebraic
number field $K$ of degree $\sigma$. So it is interesting to look
at the {\it algebraic growth characteristic} ${\mathbf
a}_K(s,r,m)$ of $f$, defined as the smallest algebraic measure of
the first $m$ derivatives of $f$ on sets $E\subset\Delta_r\cap K$
with $|E|\ge s$ (see Section \ref{S:agtf}). The following theorem,
proved in Section \ref{S:agtf}, gives lower bounds for this
characteristic.
\begin{Theorem}\label{T:tb} If $f$ is an entire function of
finite positive order then
$$\limsup_{s\to\infty}\frac{{\mathbf a}_K(s,r,m)}{s^{1/2}\log s}\ge
Cm^{1/2}.$$ If $f$ has a covering system of admissible intervals
$I(R,\alpha,\beta,\gamma,C)$, then for all $s$ sufficiently large
$${\mathbf a}_K(s,r,m)\ge
c(ms)^{1/\tau}\log\frac{ms}{a}-C_K,\;\tau=1+\frac1\gamma\,.$$
\end{Theorem}
\par Let $I_K(A)$ be the set of algebraic integers in $K$ whose
algebraic measure does not exceed $A$. Then it is possible that
$f(I_K(A))$ lies in some $I_K(B)$, like in the theorems of Polya
and Gelfond (see \cite{GS}), where $K={\Bbb Q}$ or ${\Bbb Q}[i]$.
Of course, $B$ can be large, simply due to the growth of $f$.
However, if $m(r,f)\le r^{\phi(r)}$ we prove in Section
\ref{S:agtf} the following theorem:
\begin{Theorem} If $f$ is an entire function of order $0<\rho<\sigma/2$
then
$$\liminf_{A\rightarrow\infty}
\frac{\left|I_K(A)\cap f^{-1}(I_K(\exp
A^{\phi(A)}))\right|}{|I_K(A)|}=0.$$
\end{Theorem}
This theorem tells us that, with probability close to 1, the
algebraic measure of $f(z)$ for $z\in I_K(A)$ growth faster than
$f$.
\par To prove these and other theorems, we combine our Bezout inequalities
with the standard machinery based on Siegel's lemma. This is
developed in Section \ref{S:eam} and gives lower bounds for the
algebraic measure of arguments and values of $f$ on a set
$E\subset K$.
\par We are grateful to A. Eremenko and N. Levenberg for useful
discussions. A. Eremenko also told us about Sodin's paper
\cite{So}.
\section{Characteristics of entire functions}\label{S:bf}
For an entire function $f$ on ${\Bbb C}$ we let
$$u_f(z)=\frac12\log(1+|f|^2).$$
This is a subharmonic function with Laplacian
$$\Delta u_f(z)=2\rho^2_f=\frac{2|f'|^2}{(1+|f|^2)^2}\,,$$ where
$\rho_f$ is the absolute value of the spherical derivative of $f$.
\par For a domain $D\subset{\Bbb C}$ with piecewise analytic
boundary let
\begin{equation}\begin{align}&L(D)=L_f(D)=2\int_{\partial
D}\rho_f|dz|,\notag\\
&S(D)=S_f(D)=\frac1{\pi}\int_D\rho_f^2\,d\lambda=
\frac{1}{2\pi}\int_ D\Delta u_f,\notag\end{align}\end{equation}
where $\lambda$ is the Lebesgue measure on ${\Bbb C}$. If
$D=\Delta_r$ then
$$L(D)=L(r)=2\int_{|z|=r}\rho_f|dz|\,,\,
S(D)=S(r)=\frac1\pi\int_{|z|\le r}\rho_f^2\,d\lambda$$
\par The following result of Ahlfors, with improvements by
Dufresnoy, can be found in Chapters 5 and 6 of \cite{H} and
\cite[Theorem ${\rm A}_1$, p. 190]{D}:

\begin{Theorem}\label{T:at} Let $f$ be an entire holomorphic function and let
$D$ be a domain in ${\Bbb C}$ with piecewise analytic boundary and
with Euler--Poincar\'e characteristic $\chi$. If $f$ does not
assume in $D$ the values $a\neq b$, where $|a|=|b|=1$, then
$$S(D)\le\chi+1+\frac{3}{2\delta_0}L(D),$$ where
$\delta_0$ is the spherical distance between $a$ and $b$.
Moreover, if $z\in D$ and $\operatorname{dist}(z,\partial D)=r$,
then
$$\rho_f(z)\le \frac{e^{36\pi^2/\delta_0^2}}{r}\,.$$\end{Theorem}
Here the Euler characteristic equals $-2$ for the sphere, $-1$ for
the disk, and $\chi\geq0$ for multiply connected domains.
\par The function
$$T_0(r)=T_0(r,f)=\int_0^r\frac{S(t)}t\,dt$$
is called the Ahlfors--Shimizu characteristic of $f$. If
$$m_0(r)=m_0(r,f)=\frac1{2\pi}\int_0^{2\pi}\log\sqrt{1+|f(re^{i\theta})|^2}\,d\theta,$$ then (see \cite[\S
1.5]{H})
$$T_0(r)=m_0(r)-\log\sqrt{1+|f(0)|^2}.$$ \par We let
\begin{eqnarray*}
&&M(r)=M(r,f)=\max\{|f(z)|:\,|z|=r\},\\
&&m(r)=m(r,f)=\log^+M(r,f),\\
&&T(r)=T(r,f)=\frac1{2\pi}\int_0^{2\pi}\log^+|f(re^{i\theta})|\,d\theta.\end{eqnarray*}
By \cite[Theorem 1.6]{H} and \cite[p. 13]{H} we have
\begin{eqnarray}\label{e:tm}
&&T(r)\leq m(r)\leq\frac{R+r}{R-r}\,T(R) ,\;0\leq
r<R,\\\label{e:tt0}&&\left|T(r)-T_0(r)-\log^+|f(0)|\right|\leq
\frac{\log2}{2}\,.\end{eqnarray}
\par The following relations between $L(r)$, $S(r)$ and $T_0(r)$
will be important in the sequel. If $k>1$ then
\begin{equation}\label{e:st0}
T_0(kr)\geq\int_ r^{kr} \frac{S(t)}{t}\,dt\geq S(r)\log
k.\end{equation} Moreover, H\"older's inequality implies for all
$r$
\begin{equation}\label{e:ls}L^2(r)\leq8\pi^2rS'(r).\end{equation}
\par Note that if the function $S(r)$ is bounded then $f$ is a
polynomial. Hence if $f$ is transcendental, $\epsilon>0$, $k>1$,
we can define $r_0=r_0(f,\epsilon,k)$ by
$$S(r_0)=\frac{8\pi^2}{\epsilon^2\log k}\,.$$
\begin{Lemma}\label{L:adm} If $r\geq r_0$ then there exists $r'\in(r,kr)$ so
that $L(r')\leq\epsilon S(r')$.\end{Lemma}
\begin{pf} Assuming that $L(t)>\epsilon S(t)$ for $t\in(r,kr)$, we have
by (\ref{e:ls}) $\epsilon^2S^2(t)<8\pi^2tS'(t)$. Hence
$$\frac{\epsilon^2\log k}{8\pi^2}=
\frac{\epsilon^2}{8\pi^2}\int_ r^{kr}\frac{dt}{t}< \int_
r^{kr}\frac{S'(t)}{S^2(t)}\leq\frac{1}{S(r_0)}\,,$$ a
contradiction.\end{pf}
\par We will need the following facts about functions of finite order.
Recall that (see \cite[Th.I.16]{Le}) if $\theta(r),\,r>0$, is a
positive function with
$$\rho=
\limsup_{r\to\infty}\frac{\log\theta(r)}{\log r}<\infty,$$ then
$\theta$ has a {\it proximate} order $\rho(r)$ with the following
properties: \\ (i) $\lim_{r\to\infty}\rho(r)=\rho$;\\ (ii)
$\theta(r)\le r^{\rho(r)}$, and $\theta(r_n)=r_n^{\rho(r_n)}$ for
some sequence $r_n\rightarrow\infty$;\\ (iii) the function
$\psi(r)=r^{\rho(r)-\rho}$ is {\it slowly increasing,} i.e.,
$$\lim_{r\to\infty}\frac{\psi(kr)}{\psi(r)}=1$$ uniformly on each
interval $0<a\le k\le b<\infty$. If $r^{\rho(r)-\rho}$ is a slowly
increasing function, then for every $\epsilon>0$ and every
$0<a<b<\infty$ there is $r_0$ such that
\begin{equation}\label{e:sg}(1-\epsilon)k^\rho
r^{\rho(r)}<(kr)^{\rho(kr)}<(1+\epsilon)k^\rho
r^{\rho(r)},\end{equation} for $a\leq k\leq b$ and $r\ge r_0$.
\section{Estimates for the $n$-th diameter}\label{S:end}
We will need the following lemma:
\begin{Lemma}\label{L:pms} Let $u$ be a non-negative upper bounded
subharmonic function in the disk $\Delta_R$. If $R'=R/2$ and
$\Gamma=\Delta_{R'}\cap\Delta(a,r)$, where $r<3R/4$, then
$$\frac{1}{2\pi}\,\int_{\Gamma}\Delta u\le\frac{m_0}{\log\frac{3R}{4r}}\;,$$
where $m_0=\sup\{u(z):\,|z|<R\}$. \end{Lemma}
\begin{pf} If $\Gamma\neq\emptyset$, there exists
$b\in\Delta_{R'}$ such that $\Gamma\subset\Delta(b,r)$. Since
$$0\le u(b)\le m_0+\frac{1}{2\pi}\,
\int_{\Delta_R}\log \left|\frac{R(z-b)}{R^2-\bar bz}\right|\Delta
u(z)$$ and
$$\left|\frac{R(z-b)}{R^2-\bar bz}\right|\leq
\frac{rR}{R^2-{R'}^2}=\frac{4r}{3R}$$ for $z\in\Gamma$, we obtain
$$-m_0\le\frac{\log\frac{4r}{3R}}{2\pi}\int_{\Gamma}\Delta u(z).$$
\end{pf}
\par We introduce the constants
$$\Lambda(\delta_0)=\left(4+\frac{48\pi
e^{36\pi^2/\delta_0^2}}{\delta_0}\right)^{-1},\;\Lambda=\Lambda(1).$$
The following theorem shows that the set of preimages of two
points cannot be covered by a limited number of disks of small
radius.
\begin{Theorem}\label{T:end} Let $f$ be an entire transcendental function and let
$a,b\in{\Bbb C}$, $|a|=|b|=1$, and $\delta_0$ be the spherical
distance between $a$ and $b$. If $L(R)\leq\delta_0S(R)/6$ and if
the set
$$E=\{z\in{\Bbb C}:\,|z|\leq R+r,\;z\in f^{-1}(\{a,b\})\}$$ is
covered by $n\le \Lambda(\delta_0)S(R)$ disks of radius $r$, then
$$\log \frac{3R}{4r}\leq4n\frac{m(2R)}{S(R)}\,.$$
\end{Theorem}
\begin{pf}
We assume at first that $3nr>2e^{-2}R$, so
$$\log\frac{3R}{4r}<2+\log\frac{9}{8}+\log n.$$ Using (\ref{e:st0})
and (\ref{e:tt0}) we get
\begin{equation}\label{e:stm}
S(R)\leq\frac{T_0(2R)}{\log2}\leq\frac{m(2R)}{\log2}+\frac{1}{2}\,,
\end{equation}
so
$$4n\frac{m(2R)}{S(R)}\geq4n\log2\left(1-\frac{1}{2S(R)}\right).$$
Since $n\leq\Lambda(\delta_0)S(R)\leq S(R)/4$, we obtain
$$4n\frac{m(2R)}{S(R)}\geq4n\log2\left(1-\frac{1}{8n}\right)\geq2+
\log\frac{9}{8}+\log n,$$ for all $n\geq1$. This proves Theorem
\ref{T:end} in the case $3nr>2e^{-2}R$.
\par We assume in the remainder of the proof that $3nr\le 2e^{-2}R$.
Suppose that the set $E$ can be covered by $n$ disks
$\Delta(a_j,r)$ such that $nr=d$. We claim that there are disjoint
disks $\Delta(b_j,t_j)$, $1\le j\le k\le n$, whose union contains
all disks $\Delta(a_j,2r)$ and so that $\sum t_{j}\le 2d$. For
this, we note that if the disks $\Delta(a_1,2r_1)$ and
$\Delta(a_2,2r_2)$ are not disjoint, then there is a point $z$
such that the disk $\Delta(z,2(r_1+r_2))$ contains both these
disks. Now a simple induction proves our claim.
\par We consider those disks $F_j=\Delta(b_j,t_j)$,
$j=1,\dots,l$ which intersect $\Delta_R$. Let
$\Gamma_j=F_j\cap\Delta_R$ and let
$D=\Delta_R\setminus\cup_{j=1}^l F_j$. It follows that
$$\Delta_R=D\cup\bigcup_{j=1}^l\Gamma_j,$$
thus $$S(R)=S(D)+\sum_{j=1}^lS(\Gamma_j).$$
\par By (\ref{e:stm}) and the assumption that
$S(R)\ge\Lambda^{-1}(\delta_0)$, it follows that
$m(2R)\ge\sqrt{2}$. Since $\log(1+x^2)\le4\log x$ when
$x\ge\sqrt{2}$, we get by Lemma \ref{L:pms} (with
$u=\log\sqrt{1+|f|^2}$)
$$S(\Gamma_j)\le2\frac{m(2R)}{\log \frac{3R}{2t_j}}\,.$$
Hence
$$\sum_{j=1}^lS(\Gamma_j)\le
2m(2R)\sum_{j=1}^l\frac{1}{\log \frac{3R}{2t_j}}\,.$$ Since the
sum of $2t_j/(3R)$ does not exceed $4d/(3R)\le e^{-2}$ and the
function $-1/\log x$ is concave on $(0,e^{-2})$ we conclude that
$$\sum_{j=1}^l\frac{1}{\log \frac{3R}{2t_j}}\le
\frac{l}{\log \frac{3Rl}{4d}}\;.$$ As the function $x/\log ax$ is
increasing when $x>e/a$ we have
$$\frac{l}{\log \frac{3Rl}{4d}}\le
\frac{n}{\log\frac{3R}{4r}}\,.$$ Thus
$$\sum_{j=1}^lS(\Gamma_j)\le\frac{2nm(2R)}{\log\frac{3R}{4r}}\,.$$
\par Note that the Euler characteristic $\chi_0$ of $D$
verifies $\chi_0\leq n-1$, since the domain $D$ is bounded by at
most $n+1$ Jordan curves. Moreover, we have
$$L(\partial D)\leq L(R)+\sum_{j=1}^lL(\gamma_j),$$
where $\gamma_j=\Delta_R\cap\partial\Delta(b_j,t_j)$. Thus Theorem
\ref{T:at} implies that
$$S(D)\leq n+hL(R)+h\sum_{j=1}^lL(\gamma_j),$$where $h=3/(2\delta_0)$.
If $z\in\gamma_j$, then $f$ does not take the values $a$ and $b$
in the disc $\Delta(z,r)$, so by Theorem \ref{T:at} $\rho_f(z)\leq
h_1/r$, where $h_1=e^{36\pi^2/\delta_0^2}$. Hence
$L(\gamma_j)\leq4\pi h_1 t_j/r$ and
$$\sum_{j=1}^lL(\gamma_j)\le\frac{8\pi h_1d}{r}=8\pi h_1n.$$
We conclude that $S(D)\leq (1+8\pi hh_1)n+hL(R)$, so
$$S(R)\le (1+8\pi hh_1)n+hL(R)+\frac{2nm(2R)}{\log\frac{3R}{4r}}\,.$$ If
$n\leq\Lambda(\delta_0)S(R)$ and $L(R)\leq\delta_0S(R)/6$, then
$$\log\frac{3R}{4r}\le \frac{4nm(2R)}
{S(R)}\,.$$
\end{pf}
\par For a set $G\subset{\Bbb C}$ and an integer $n\geq1$
we introduced in \cite{CP2} the {\it $n$-th diameter } of $G$ as
$$\operatorname{diam}_n(G)=\inf\left\{r_1+\dots+r_k:\,k\leq
n,\;G\subset\bigcup_{j=1}^kC_j(r_j)\right\},$$ where $C_j(r_j)$
are closed disks of radii $r_j>0$.
\par Given a non-constant entire function $f$ we denote in the sequel
by $n_0=n_0(f)$ the maximum of the numbers
$|f^{-1}(w)\cap\Delta_2|$ when $|w|=1$.
\begin{Corollary}\label{C:end} In the assumptions of Theorem \ref{T:end},
let
$$F=\{z\in{\Bbb C}:\,2\leq|z|\leq
R+1,\;z\in f^{-1}(\{a,b\})\}.$$ If $L(R)\leq\delta_0S(R)/6$,
$n\leq\Lambda(\delta_0)S(R)-2n_0$, and
$d_n=\operatorname{diam}_n(F)<1$, then
$$\log\frac{R}{d_n}\le 4(n+2n_0)\frac{m(2R)}{S(R)}+\log\frac{4}{3}\,.$$
\end{Corollary}
\begin{pf} If $\epsilon>0$ and $d_n+\epsilon<1$, we can cover $F$
by $n$ disks of radius $d_n+\epsilon$. The number of points of
$f^{-1}(\{a,b\})\cap\Delta_2$ does not exceed $2n_0$. We cover
them with $2n_0$ disks of radius $d_n+\epsilon$. Since
$d_n+\epsilon<1$, we apply Theorem \ref{T:end} and then let
$\epsilon\rightarrow0$.
\end{pf}
\par Let $$D_n(\theta,r)=\{z\in{\Bbb C}:\,2\leq|z|\leq
r,\;f(z)=e^{i\theta}\}$$ and
$$d_n(\theta,r)=\min\{1,{\rm diam}_n(D_n(\theta,r))\}.$$
\begin{Corollary}\label{C:diam}
Let $f$ be an entire transcendental function. If $L(R)\leq S(R)/6$
and $n\le\frac{1}{2}\,\Lambda S(R)-n_0$, then
$$\log\frac{R}{d_n(\theta,R+1)}\leq
\max\left\{8(n+n_0)\frac{m(2R)}{S(R)} +\log3,\,\log(2R)\right\}$$
for all $e^{i\theta}$ in an arc of length $l>\pi$ in
$\partial\Delta$.
\end{Corollary}
\begin{pf} Suppose that the spherical distance between $a=e^{i\phi}$
and $b=e^{i\psi}$ is at least $\delta_0=1$. Let $m=2n$. If $F$ is
as in Corollary \ref{C:end}, then
$$\operatorname{diam}_{m}(F)\le d_{n}(\phi,R+1)+d_{n}(\psi,R+1).$$ Hence if
$\operatorname{diam}_{m}(F)<1$ we have by Corollary \ref{C:end}
$$-\log\left(\frac{d_n(\phi,R+1)}{R}+
\frac{d_n(\psi,R+1)}{R}\right)\leq
4(m+2n_0)\frac{m(2R)}{S(R)}+\log\frac{4}{3}\,.$$ If
$0<\alpha\le\beta$, then $\log(\alpha+\beta)\le\log\beta+\log2$.
Thus
$$\log\frac{R}{\max\{d_n(\phi,R+1),d_n(\psi,R+1)\}}\leq
8(n+n_0)\frac{m(2R)}{S(R)}+\log3.$$ If
$\operatorname{diam}_{m}(F)\geq1$ then
$$\max\{d_n(\phi,R+1),d_n(\psi,R+1)\}\geq1/2.$$
\par Consequently, if the estimate in the statement of the corollary
fails for some $e^{i\phi}$, then it must hold for all $e^{i\psi}$
lying at spherical distance at least 1 from $e^{i\phi}$. Since the
set of such $e^{i\psi}$ is an arc of length greater than $\pi$,
the corollary follows.
\end{pf}
\section{General estimates for $e_n$ and $m_n(r)$}\label{S:gpe}
Let $f$ be an entire transcendental function and recall that
$n_0=n_0(f)$ is the maximum of the numbers
$|f^{-1}(w)\cap\Delta_2|$, when $|w|=1$. In the following lemma,
the estimates on the $n$-th diameter obtained in the previous
section, combined with results form \cite{CP2}, lead to estimates
of the transcendence measures $e_n$ and $m_n(r)$ in terms of
$m(r)$.
\begin{Lemma}\label{L:pe} Let $R_0=R_0(f)$ be the largest
among the unique solutions of the equations:
$$R=64,\;\;\;S(R)=\frac{288\pi^2}{\log(4/3)}\,,\;\;\;
m(R)=4\log^+R,\;\;\;m(4R)=36.$$ If $R>R_0$ and
$n\le\frac{1}{2}\,\Lambda S(R)-n_0$, then
\begin{eqnarray*}e_n&\leq&2nm(4R)\log R,\\
m_n(r)&\leq&3nm(4R)\log r, 1\leq r\leq R.\end{eqnarray*}
\end{Lemma}
\begin{pf} Using Lemma \ref{L:adm} with $k=4/3$ and $\epsilon=1/6$
we can find, for all $R>R_0$, a radius $R'\in(R,4R/3)$ so that
$L(R')\leq S(R')/6$. Since $n+n_0\le\frac12\Lambda S(R)$, we have
by Corollary \ref{C:diam} with $r=R'+1$, that for $e^{i\theta}$ in
a set of length $l>\pi$ in $\partial\Delta$
$$\log\frac{R'}{d_n(\theta,r)}\leq
\max\left\{4\Lambda m(2r)+\log3,\,\log(2r)\right\}.$$ Since
$er<3R'$, the latter inequality implies
$$\log\frac{36er}{d_n(\theta,r)}\leq5+\max\left\{4\Lambda m(2r)+\log
3,\,\log(2r)\right\}.$$
\par Theorem 4.2 in \cite{CP2} asserts that if for some $r\geq2$
one has $d_n(\theta,r)\geq a$ on a set $E\subset\partial\Delta$ of
length $l$, then $$e_n\leq n\max\{m(er),\log(er)\}\log r+
n\log(er)\left(\log\frac{36er}{a}+\frac{4\pi}{l}\right).$$
\par Suppose that $4\Lambda m(2r)+\log3\geq\log(2r)$.
Since $m(er)\geq\log(er)$, the above estimate yields
$$e_n\le nm(er)\log r+n\log(er)(11+4\Lambda m(2r)).$$
Since $er<4R$, $R>64$, $m(4R)>36$ and $\Lambda<e^{-300}$ we have
\begin{eqnarray*}
e_n&\leq&\left(1+4\Lambda+\frac{11}{m(4R)}\right)nm(4R)\log(4R)\\
&\leq&\frac{4}{3}\left(1+4\Lambda+\frac{11}{36}\right)nm(4R)\log
R<2nm(4R)\log R.
\end{eqnarray*}
\par If $4\Lambda m(2r)+\log3<\log(2r)$, then using in addition that $4\log(4R)\leq m(4R)$,
we get
\begin{eqnarray*}
e_n&\leq&nm(er)\log r+n\log(er)(9+\log(2r))\\
&\leq&\left(1+\frac{9}{m(4R)}+\frac{1}{4}\right)nm(4R)\log(4R)\leq2nm(4R)\log
R.
\end{eqnarray*}
\par Thus $e_n\leq2nm(4R)\log R$. By \cite[\S 4
(5)]{CP2} we have for $1\leq r\leq R$
$$m_n(r)\le\frac{e_n+nm(R)}{\log R}\,\log r\leq3nm(4R)\log r.$$
\end{pf}
\par The above lemma shows that estimates for $e_n$ and $m_n$
require knowledge of the relationship between $m(4R)$ and $S(R)$.
The following theorem shows the kind of hypotheses on $m(4R)$ and
$S(R)$ needed to get good estimates on $e_n$ and $m_n$.
\par We denote by $R_1(f)$ the maximum of $R_0(f)$ and the solution of the
equation
$$T_0(r)=(3\log2)/2+3\log^+|f(0)|.$$We call an interval
$$I(R,\alpha,\beta,\gamma,C)=\left[\beta S^\gamma(R),\frac12\,\Lambda S(R)-n_0\right]$$ {\it
admissible} if $R>R_1(f)$, $\alpha,\beta>0$, $0<\gamma\le1$,
$\beta S^\gamma(R)\le\frac12\,\Lambda S(R)-n_0-1$, $S(R)\ge
R^\alpha$ and $m(4R)\le CS(R)$. We let $I(\alpha,\beta,\gamma,C)$
be the union of all admissible intervals
$I(R,\alpha,\beta,\gamma,C)$.
\begin{Theorem}\label{T:me1} If $n\in I(\alpha,\beta,\gamma,C)$, then
$$e_n\leq \frac{2C}{\alpha\gamma\beta^{1/\gamma}}\,n^{1+1/\gamma}
\log\frac{n}{\beta}\,.$$ If $R$ is so that $n\in
I(R,\alpha,\beta,\gamma,C)$, then $2\Lambda m(2R)\ge n$ and
$$m_n(r)\leq\frac{3C}{\beta^{1/\gamma}}\, n^{1+1/\gamma}\log
r,\;1\leq r\leq R.$$\end{Theorem}
\begin{pf} By the properties of admissible intervals we have
$$m(4R)\le CS(R)\le C\left(\frac n\beta\right)^{1/\gamma}$$
and $$R\le S^{1/\alpha}(R)\le\left(\frac
n\beta\right)^{1/(\alpha\gamma)}.$$ By Lemma \ref{L:pe}
$$e_n\leq
\frac{2C}{\alpha\gamma\beta^{1/\gamma}}\,n^{1+1/\gamma}
\log\frac{n}{\beta}$$ and
$$m_n(r)\leq\frac{3C}{\beta^{1/\gamma}}\,n^{1+1/\gamma}\log
r,\;1\leq r\leq R.$$
\par By (\ref{e:tm}) and (\ref{e:tt0}) we have
$$T_0(r)-\frac{\log2}{2}\leq m(r)\leq
3T_0(2r)+\frac{3\log2}{2}+3\log^+|f(0)|.$$ So if $r\geq R_1(f)$
then
\begin{equation}\label{e:mt0}
\frac{1}{2}\,T_0(r)\le m(r)\le4T_0(2r).
\end{equation}
By (\ref{e:st0}) $S(R)\log 2\le T_0(2R)$, so
\begin{equation}\label{e:sm}S(R)\leq\frac{2m(2R)}{\log 2}\,.\end{equation}
Therefore $n\leq\Lambda S(R)/2\le 2\Lambda m(2R)$, and the proof
is complete.\end{pf}
\par The following corollary establishes a case when polynomial
estimates for $e_n$ hold for all $n$.
\begin{Corollary}\label{C:me2} Suppose that there is a sequence of admissible
intervals $I(R_j,\alpha,\beta,\gamma,C)$ such that
$R_j\rightarrow\infty$ and $\beta S^\gamma(R_{j+1})\le\Lambda
S(R_j)/2-n_0$, $j\ge1$. Then the conclusions of Theorem
\ref{T:me1} hold for all $n\ge\beta S^\gamma(R_1)$.
\end{Corollary}
\par A system of admissible intervals satisfying the hypotheses of
this corollary will be called a {\em covering system}.
\section{The lower order of transcendence}\label{S:loot}
\par In order to apply Theorem \ref{T:me1} effectively we need
information on the set $I(\alpha,\beta,\gamma,C)$. The main goal
of this section is to establish that for every entire function of
finite positive order we can find $\alpha$, $\beta$ and $C$ such
that the set $I(\alpha,\beta,1,C)$ is unbounded. Then Theorem
\ref{T:me1} will imply that the lower order of transcendence
$\underline\tau(f)$ of any entire function $f$ of finite positive
order is 2. Our first step is to study the ratio of $T_0(r)$ and
$S(r)$.
\begin{Lemma}\label{L:mi} If $AT_0(r_1/k)\leq T_0(r_1)$, where
$k>1$ and $A>1$, then there is $r\in(r_1/k,r_1)$ such that
$cS(r)\ge T_0(r_1)$, where
$$c=\frac{A\log k}{\log A}\,.$$ \end{Lemma}
\begin{pf} Let us take $r_2$ such that $T_0(r_1)=AT_0(r_2)$.
Then $r_1/k\le r_2<r_1$. If
$$T_0(t)>\frac{\log k}{\log A}\,S(t)$$
on $(r_2,r_1)$, then
$$\log A=\int_{r_2}^{r_1}\frac{T_0'(t)}{T_0(t)}\,dt=
\int_{r_2}^{r_1}\frac{S(t)}{tT_0(t)}\,dt< \frac{\log A}{\log
k}\,\int_{r_2}^{r_1}\frac{dt}{t}\leq\log A.$$ Hence there is
$r\in(r_2,r_1)$ such that $$\frac{\log k}{\log A}\,S(r)\ge
T_0(r)\ge T_0(r_2)=\frac{T_0(r_1)}A\,.$$
\end{pf}
\par Next we need the ratio $m(4r)/S(r)$ to be bounded above for some
numbers $r$. The following lemma provides sufficient conditions
for such values of $r$.
\begin{Lemma}\label{L:pl} Suppose that for some $k>1$ and $r_1>kR_1(f)$ there
are constants $A_1>8$ and $A_2>1$ such that $A_1m(r_1/k)\le
m(r_1)$ and $m(8r_1)\le A_2m(r_1)$. Then there is
$r\in(r_1/k,2r_1)$ such that $CS(r)\ge m(8r_1)$, where
$$C=\frac{A_1A_2\log(2k)}{2\log(A_1/8)}\,.$$
\end{Lemma}
\begin{pf} By (\ref{e:mt0}) it follows that
$$\frac{1}{2}\,T_0(r_1/k)\le
m(r_1/k)\le\frac{m(r_1)}{A_1} \le\frac{4T_0(2r_1)}{A_1}\,.$$ By
Lemma \ref{L:mi} there is $r\in(r_1/k,2r_1)$, such that
$c_1S(r)\ge T_0(2r_1)$, where $$c_1=\frac{A_1\log(2k)}{8\log
(A_1/8)}\,.$$ Hence
$$c_1S(r)\ge T_0(2r_1)\ge\frac{m(r_1)}4\ge\frac{m(8r_1)}{4A_2}\,.$$
\end{pf}
\par As we will now see, the ratio $m(4r)/S(r)$ is bounded
near points where $m(r)$ is close to its proximate order.
\begin{Lemma}\label{L:si} Suppose that $m(r)\le r^{\phi(r)}$, where
$\lim_{r\to\infty}\phi(r)=\rho$ and the function
$r^{\phi(r)-\rho}$ is slowly increasing. Let $0<a\le1$ and
$k^\rho>16/a$. If $r_1$ is sufficiently large and $m(r_1)\ge
ar_1^{\phi(r_1)}$, then there is $r\in(r_1/k,2r_1)$ such that
$CS(r)\ge m(8r_1)$, where
$$C=\frac{(8k)^\rho\log(2k)}{2\log (ak^\rho/16)}\,.$$
\end{Lemma}
\begin{pf}``Sufficiently large" in the statement of the lemma means
$r_1>kr_0$, where $r_0>R_1(f)$ is a number such that
$$\frac{1}{2}\,b^\rho r^{\phi(r)}<(br)^{\phi(br)}<2b^\rho
r^{\phi(r)}$$ holds for $k^{-1}\leq b\leq 8$ and $r\ge r_0$ (see
(\ref{e:sg})). Then
$$m(r_1/k)\le\left(\frac{r_1}k\right)^{\phi(r_1/k)}\le
2k^{-\rho}r_1^{\phi(r_1)}\le \frac{2}{ak^\rho}\,m(r_1),$$ and
$$m(8r_1)\le(8r_1)^{\phi(8r_1)}\le
2^{3\rho+1}r_1^{\phi(r_1)}\le\frac{2^{3\rho+1}}{a}\,m(r_1).$$ The
conclusion follows by Lemma \ref{L:pl}, if we let $A_1=ak^\rho/2$
and $A_2=2^{3\rho+1}/a$.
\end{pf}
\par In the following theorem we prove that $\underline\tau(f)=2$. Note that
we also give effective estimates on the ``type" of growth of $e_n$
and $m_n(r)$.
\begin{Theorem}\label{T:liminf} Let $f$ be an entire function of finite
order $\rho>0$. There exist sequences of integers
$n_j\nearrow\infty$ and $\epsilon_j\rightarrow0,\;\epsilon_j>0$,
such that
\begin{eqnarray*}
&&\frac{n_j^2\log n_j}{2\rho+1}\leq e_{n_j}\leq
\frac{8^{\rho+3}(\rho+5)}{\Lambda\rho^2}\,n_j^2\log n_j\,,\\
&&m_{n_j}(r)\leq\frac{8^{\rho+3}(\rho+5)}{\Lambda\rho}\,n_j^2\log
r, \;1\leq r\leq\frac{1}{2}\,n_j^{1/\rho-\epsilon_j}\,.
\end{eqnarray*}\end{Theorem}
\begin{pf} Let $\rho(r)$ be a proximate order for $m(r)$. By
its definition there exists a sequence $R'_j\rightarrow\infty$
such that $m(R'_j)=(R'_j)^{\rho(R'_j)}$. Take $a=1$ and
$k=2^{5/\rho}$. By Lemma \ref{L:si} there exist, for all $j$
sufficiently large, numbers $R_j\in(R'_j/k,2R'_j)$ such that
$$CS(R_j)\ge m(8R'_j)\ge m(4R_j),$$
where
\begin{equation}\label{e:c}
C=\frac{2^{3\rho+4}(\rho+5)}{\rho}\,.\end{equation} Since for $j$
large
$$m(8R'_j)>(R'_j)^{\rho(R'_j)}>(R_j/2)^{3\rho/4},$$
we see that  $S(R_j)\ge R_j^{\rho/2}$. Also $S(R_j)\ge
6(n_0+1)/\Lambda$ when $j$ is sufficiently large. \par Hence the
intervals $I_j=I(R_j,\rho/2,\Lambda/3,1,C)$ are admissible and
there exists a sequence of integers $n_j\in I_j$. By Theorem
\ref{T:me1}
$$e_{n_j}\leq\frac{3\cdot2^{3\rho+6}(\rho+5)}
{\Lambda\rho^2}\,n_j^2\log\frac{3n_j}{\Lambda}\leq
\frac{2^{3\rho+9}(\rho+5)}{\Lambda\rho^2}\,n_j^2\log n_j,$$ for
all $j$ sufficiently large. Moreover
$$m_{n_j}(r)\leq\frac{2^{3\rho+8}(\rho+5)}{\Lambda\rho}\,n_j^2 \log
r,\;1\leq r\leq R_j,$$ and
$$n_j\le 2\Lambda m(2R_j)\le2\Lambda(2R_j)^{\rho(2R_j)}.$$ By (\ref{e:sg})
there is a sequence of positive $\epsilon_j\to0$ such that
$$n_j\le2\Lambda(2R_j)^{\rho(2R_j)}\le
4\Lambda2^\rho R_j^{\rho(R_j)}\le2^\rho
R_j^{1/(1/\rho-\epsilon_j)}.$$ Hence $R_j\ge
n_j^{1/\rho-\epsilon_j}/2$.
\par For the lower estimate on $e_{n_j}$, we take $r$ with
$n_j=m(r)<r^{\rho+1/4}$, so $\log r>\log n_j/(\rho+1/4)$. By
\cite[\S 4 (3)]{CP2} and \cite[Corollary 2.6]{CP2}
$$e_{n_j}\geq\frac{n_j^2}{2}\,\log r-n_jm(r)\geq \frac{n_j^2\log
n_j}{2\rho+1/2}-n_j^2\geq\frac{n_j^2\log n_j}{2\rho+1}\,.$$
\end{pf}
\section{Doubling inequalities}\label{S:di}
In this section we prove {\it doubling inequalities}, which provide
upper bounds for the ratio $M(2r,F)/M(r,F)$, where
$F(z)=P(z,f(z))$ and $P$ is a polynomial. For $f(z)=e^z$
such inequalities were obtained by Tijdeman in \cite{Ti1}.
\par We will need a simple lemma, whose proof is
contained in the proof of Theorem 2.2 of \cite{CP2}.
\begin{Lemma}\label{L:bp} If $r<s$ and an entire function $f$ has $m$
zeros in $\Delta_r$, then
$$\frac{M(s,f)}{M(r,f)}\ge\left(\frac{r^2+s^2}{2rs}\right)^m.$$
\end{Lemma}
\par First, we reduce the problem of doubling inequalities to
the problem of obtaining estimates for the transcendence measures
$m_n$ of dilations $f_r(z)=f(rz)$ of $f$. Let $r\ge 1$, let
$P(z,w)$ be a polynomial of degree at most $n$ and let
$F(z)=P(z,f(z))$. If $P_r(z,w)=P(rz,w)$, then
$\|P_r(z,f_r(z))\|_\Delta=M(r,F)$, while
$\|P_r(z,f_r(z))\|_{\Delta_2}=M(2r,F)$. Hence
$$\frac{M(2r,F)}{M(r,F)}\le e^{m_n(2,f_r)},$$
and we have to estimate $m_n(2,f_r)$. \begin{Theorem}\label{T:gdi}
Let $f$ be an entire transcendental function of finite positive
order $\rho$. There exists a sequence of integers $\{n_j\}$
increasing to infinity with the following property: For every
$r\ge1$ there is an integer $j_r$ such that
$$\frac{M(2r,F)}{M(r,F)}\le 2^{an_j^2},\;\;\;
a=\frac{8^{\rho+3}(\rho+5)}{\Lambda\rho}\,,\;\;\; j\ge j_r,$$
where $F(z)=P(z,f(z))$ and $P(z,w)$ is a polynomial of degree at
most $n_j$.
\end{Theorem}
\begin{pf}
\par Let us denote by $n_r$ the maximum of the
numbers $|f_r^{-1}(w)\cap\Delta_2|$, when $|w|=1$. Let $w_r$ be a
point where this maximum is achieved and let
$g_r(z)=f(rz)-w_r$. By Lemma \ref{L:bp}
$$\left(\frac{5}{4}\right)^{n_r}\leq\frac{M(4,g_r)}{M(2,g_r)}\,.$$ Since $f$ is
not constant, there exists $\epsilon>0$ such that $M(2,g_r)\geq\epsilon$
for every $r\geq1$. Since $M(4,g_r)\leq M(4r)+1$ it follows that
\begin{equation}\label{e:nr}
n_r\leq C_1m(4r)-1,\end{equation} where $C_1$ is a constant depending
only on $f$.
\par Let $I=I(R,\alpha,\beta,\gamma,C)$ be an admissible interval for
$f$. From the definition of the number $R_1(f)$ in Section
\ref{S:gpe} it follows that $R_1(f_r)\le R_1(f)$. Note that
$m(t,f_r)=m(rt)$ and $S_{f_r}(t)=S(rt)$. Hence, if $R'=R/r$, then
$S_{f_r}(R')\ge{R'}^\alpha$ and $m(4R',f_r)\le CS_{f_r}(R')$.
Therefore, the interval $I'=(R',\alpha,\beta,\gamma,C)$ is
admissible for $f_r$ if $R'\ge R_1(f)$ and
$$\frac{\Lambda}2\,S(R)-n_r-1\ge \beta S^\gamma(R).$$
\par Since $f$ has finite positive order $\rho$, by the proof of
Theorem \ref{T:liminf} there is a sequence $R_j$ increasing to infinity
such that the intervals
$$I_j=I(R_j,\rho/2,\Lambda/3,1,C)=
\left[\frac{\Lambda S(R_j)}{3}\,,\frac{\Lambda S(R_j)}{2}-n_0\right]$$
are admissible, where $C$ is defined in (\ref{e:c}). For every $r\ge 1$
let $j_r$ be the smallest integer such that $R_{j_r}>rR_1(f)$ and
$$\frac{\Lambda}{10}\,S(R_{j_r})\geq C_1m(4r).$$
Then
$$\frac{\Lambda}2\,S(R_{j_r})-n_r-1\ge\frac{2\Lambda}{5}\,S(R_{j_r})\ge
\frac\Lambda3\,S(R_{j_r})$$ and the intervals
$$I'_j=I(R_j/r,\rho/2,\Lambda/3,1,C)=
\left[\frac{\Lambda S(R_j)}{3}\,,\frac{\Lambda S(R_j)}{2}-n_r\right]$$
are admissible for $f_r$ when $j\ge j_r$.
\par Let $j_0$ be the smallest integer so that
$S(R_{j_0})\geq\max\{15/\Lambda,10n_0/\Lambda\}$. Then for $j\geq
j_0$ the intervals $I_j''=[\Lambda S(R_j)/3,2\Lambda S(R_j)/5]$
contain an integer $n_j$ and $I_j''\subset I_j$. Moreover, if
$j\geq j_r$ then $I''_j\subset I'_j$, so by Theorem \ref{T:me1}
$$m_{n_j}(2,f_r)\le\frac{9C}\Lambda\,n_j^2\log 2\leq an_j^2\log 2.$$
Consequently $M(2r,F)/M(r,F)\le 2^{an_j^2}$, for all $j\ge
j_r$.\end{pf} \vspace{2mm}\noindent{\bf Remark.} With the
notations of the above proof, since $n_j\in I_j''\subset I_j$ it
follows that the conclusions of Theorem \ref{T:liminf} hold for
the sequence $\{n_j\}$ constructed in Theorem \ref{T:gdi}. A
sequence of integers $\{n_j\}$ increasing to infinity for which
the conclusions of both Theorems \ref{T:liminf} and \ref{T:gdi}
are valid, will be called a {\em fundamental sequence} for
$f$.\vspace{2mm}
\par In the following theorem we prove doubling inequalities
for functions which possess a covering system of admissible
intervals.
\begin{Theorem}\label{T:di}Let $f$ be an entire transcendental function
which has a covering system of admissible intervals
$I_j=I(R_j,\alpha,\beta,\gamma,C)$. For every $r\ge1$ there exists
an integer $j_r$ such that
$$\frac{M(2r,F)}{M(r,F)}\leq\left\{\begin{array}{ll}
\exp\left(3nm(4R_{j_r})\log 2\right),\;\mbox{if $n<\beta
S^\gamma(R_{j_r})/2$},\\ \\
\exp\left(3C(2\beta^{-1})^{1/\gamma}n^{1+1/\gamma}\log
2\right),\:\mbox{if $n\ge\beta S^\gamma(R_{j_r})/2$},\end{array}
\right.$$ where $F(z)=P(z,f(z))$ and $P(z,w)$ is a polynomial of
degree at most $n$.\end{Theorem}
\begin{pf} Let $j_r$ be the smallest integer such that
\begin{equation}\label{e:jr}
R_{j_r}>rR_1(f),\;\Lambda m(4R_{j_r})\ge 4CC_1m(4r),
\end{equation} where $C_1$ is the constant from (\ref{e:nr}).
By (\ref{e:nr}) and the properties of admissible intervals we have for
$j\geq j_r$
$$\frac{\Lambda}{4}\,S(R_j)\geq\frac{\Lambda}{4C}\,m(4R_{j_r})
\geq C_1m(4r)\geq n_r+1,$$
so $$\frac{\Lambda}{2}\,S(R_j)-n_r-1\geq
\frac{\Lambda}{4}\,S(R_j)\geq \frac{\beta}{2}\,S^\gamma(R_j).$$
Moreover, since $I_j$ form a covering system we have
$$\frac\beta2\,S^\gamma(R_{j+1})\leq\frac\Lambda4\,
S(R_j)-\frac{n_0}{2}\leq\frac\Lambda2\,S(R_j)-n_r.$$ Thus the
intervals $I'_j=I(R_j/r,\alpha,\beta/2,\gamma,C)$, $j\geq j_r$,
form a covering system of admissible intervals for $f_r$. By
Corollary \ref{C:me2}
$$m_n(2,f_r)\leq3C(2\beta^{-1})^{1/\gamma}\, n^{1+1/\gamma}\log
2,$$when $n\ge\beta S^\gamma(R_{j_r})/2$.
\par If $$n<\beta S^\gamma(R_{j_r})/2\le\frac\Lambda2\,S(R_{j_r})-n_r,$$
then by Lemma \ref{L:pe} $m_n(2,f_r)\le 3nm(4R_{j_r})\log2$.
\end{pf}
\par Let us denote by $Z_n(r,f)=Z_n(r)$ the maximal number
of zeros of $P(z,f(z))$ in $\Delta_r$, when $P(z,w)$ is a
polynomial of degree at most $n$. In Corollary 2.6 of \cite{CP2}
we proved that $Z_n(r)\le 2m_n(3r)$. Now we can improve this
estimate. \par The first result gives an estimate on $Z_n(r)$ for
all transcendental functions of finite positive order. Note that
the constant $a$ depends only on the order $\rho$ of $f$.
\begin{Corollary}\label{C:gnz} If $\{n_j\}$ is a fundamental sequence for
$f$ then $Z_{n_j}(r)\le 4an_j^2$, for $r\geq1$ and $j\ge j_r$.
\end{Corollary}
\begin{pf} Let $P(z,w)$ be a polynomial of degree $n_j$ such that the
number of zeros of $F(z)=P(z,f(z))$ in $\Delta_r$ equals
$Z_{n_j}(r)$. Then by Theorem \ref{T:gdi} and Lemma \ref{L:bp}
$$\left(\frac54\right)^{Z_{n_j}(r)}\le\frac{M(2r,F)}{M(r,F)}\le2^{an_j^2}$$
when $j\ge j_r$. Hence
$$Z_{n_j}(r)\le\frac{an_j^2\log2}{\log(5/4)}\le4an_j^2.$$
\end{pf}
\par The second corollary provides estimates on $Z_n(r)$ for all $n$
and has a similar proof. \begin{Corollary}\label{C:nz} In the
assumptions of Theorem \ref{T:di} we have
$$Z_n(r)\leq\left\{\begin{array}{ll}10nm(4R_{j_r}),\;
\mbox{if $n<\beta S^\gamma(R_{j_r})/2$},\\ \\
10C(2\beta^{-1})^{1/\gamma}\,n^{1+1/\gamma},\; \mbox{if
$n\geq\beta
S^\gamma(R_{j_r})/2$}.\end{array}\right.$$\end{Corollary}
\par Doubling inequalities lead to tangential Markov inequalities,
which provide upper estimates for the derivative of the function
$F(z)=P(z,f(z))$, where $P(z,w)$ is a polynomial of degree $n$. As
before, we give two versions of such inequalities: one for general
entire functions and another for functions with a covering system
of admissible intervals.
\begin{Theorem}\label{T:gmi} Let $\{n_j\}$ be a fundamental sequence for $f$.
For every $r\ge1$ there is an integer $j_r$ such that
$$M(r,F')\le \frac{eaM(r,F)n_j^2}r\,,$$ where $F(z)=P(z,f(z))$, $P(z,w)$
is a polynomial of degree $n_j$, and $j\ge j_r$.
\end{Theorem}
\begin{pf}\par For $r\ge1$ let $j_r$ be the integer from Theorem
\ref{T:gdi}. Let $F(z)=P(z,f(z))$, where $P(z,w)$ is a polynomial
of degree $n_j$ and $j\ge j_r$. Since $m(r,F)$ is a convex increasing
function of $\log r$
$$|F(z)|\leq M(r,F)\exp\left((m(t,F)-m(r,F))\frac{\log(|z|/r)}{\log
(t/r)}\right), \ r\le|z|\le t.$$ Let $b>1$ be such that
$$\frac{m(t,F)-m(r,F)}{\log (t/r)}\leq b.$$ The function
$h(x)=e^{b\log(1+x)}/x$ attains its minimum value when $x=x_b=1/(b-1)$,
and $h(x_b)<eb$. Therefore, if $r(1+x_b)\le t$ and
$|z|=r$, then by the Cauchy estimates
$$|F'(z)|\le \frac{M(r,F)}{rx_b}\,e^{b\log(1+x_b)}\le
\frac{ebM(r,F)}{r}\,.$$ Taking $t=2r$, we have by Theorem
\ref{T:gdi} $(m(2r,F)-m(r,F))/\log2\le b=an_j^2$ and $1+x_b\leq2$. Thus
$$|F'(z)|\le \frac{eaM(r,F)n_j^2}r\,.$$
\end{pf}
\par The following theorem provides estimates on $M(r,F')$ for all $n$
and has a similar proof.
\begin{Theorem}\label{T:mi} In the assumptions of Theorem \ref{T:di} we have
$$M(r,F')\le \frac{3eC2^{1/\gamma}n^{1+1/\gamma}M(r,F)}{\beta^{1/\gamma}r}\,,$$
where $F(z)=P(z,f(z))$, $P(z,w)$ is a polynomial of degree $n$,
and $n\geq\beta S^\gamma(R_{j_r})/2$.\end{Theorem}
\section{Special classes of functions}\label{S:scof}
\par In this section we find sufficient conditions
for a function $f$ to have estimates of the form $e_n=O(n^\tau\log
n)$ for some $\tau\geq2$. These conditions are imposed on the
growth of $f$ and are easy to verify. We start with the class of
entire functions $f$ whose growth satisfies the following
inequalities: There exist constants $A_2>A_1>1$ and $k>1$ such
that
\begin{equation}\label{e:sc1} A_1m(r)\le m(kr)\le A_2m(r)
\end{equation}
for all $r$ sufficiently large. These are functions of finite
positive order and this class includes, for example, all functions
$$f(z)=\sum_{j=1}^mp_j(z)e^{q_j(z)},$$ where $p_j$ and $q_j$ are
polynomials. Moreover we show in the next section that the Riemann
$\zeta$-function and the function $\xi$ are also in this class.
\begin{Theorem}\label{T:sc1} Let $f$ be an entire function of
order $\rho$ which satisfies (\ref{e:sc1}) for all $r$
sufficiently large. Then, for all $n$ sufficiently large,
$$e_n\le K_1n^2\log n\;,
\;m_n(r)\le K_2n^2\log r,\;1\le r\le n^{1/\rho-\epsilon_n}/2,$$
where the constants $K_1,K_2$ depend only on $A_1,A_2,k,$ and
$\epsilon_n>0,\epsilon_n\to0$.
\end{Theorem}
\begin{pf} Inequalities (\ref{e:sc1}) imply that
$$A_1^jm(r)\leq m(k^jr)\leq A_2^jm(r),\;
d_1r^{\rho_1}\le m(r)\le d_2r^{\rho_2},$$ where
$$\rho_1=\frac{\log A_1}{\log k}\,,\;\rho_2=\frac{\log A_2}{\log k}\,,\;
d_1=\frac{m(1)}{A_1}\,,\;d_2=A_2m(1).$$ Thus $f$ is a function of
finite positive order $\rho\in[\rho_1,\rho_2]$, and we may assume
that (\ref{e:sc1}) holds with constants $k\geq8$ and $A_1>8$. Then
$A_1m(r/k)\le m(r)$ and $m(8r)\le A_2m(r)$.
\par For every $r$ sufficiently large there is, by Lemma \ref{L:pl},
$r'\in(r/k,2r)$ such that $CS(r')\ge m(4r')$, where
$$C=\frac{A_1A_2\log(2k)}{2\log (A_1/8)}\,.$$
In particular, for all $j$ sufficiently large, there is
$R_j\in((2k)^{j},(2k)^{j+1})$ such that $CS(R_j)\ge m(4R_j)$.
Since the order of $f$ is $\rho$ we may assume that $S(R_j)\ge
R_j^{\rho/2}$.
\par Using (\ref{e:sm}) and (\ref{e:sc1}) we get
$$S(R_{j+1})\leq3m(2R_{j+1})\leq
3m(8k^2(2k)^j)\leq3A_2^3m(R_j)\leq3CA_2^3S(R_j).$$ Hence
$S(R_{j+1})/S(R_j)\leq M=3CA_2^3$. \par Let $j$ be so large that
$S(R_j)\geq6(n_0+1)/\Lambda$ and let $\beta=\Lambda/(3M)$. Then
$$\beta S(R_j)\leq\beta S(R_{j+1})\leq\frac{\Lambda}{3}\,
S(R_j)\leq\frac{\Lambda}{2}\,S(R_j)-n_0-1,$$ so the intervals
$I_j=I(R_j,\rho/2,\beta,1,C)$ form a covering system of admissible
intervals, starting with some $j$ sufficiently large. The theorem
now follows from Corollary \ref{C:me2}. If $n\in I_{j_n}$ then the
fact that $R_{j_n}\ge n^{1/\rho-\epsilon_n}/2$ can be proved
exactly like the similar statement in Theorem \ref{T:liminf}.
\end{pf}
\par The functions $f$ satisfying (\ref{e:sc1}) have covering
systems of admissible intervals. Hence they also satisfy the
hypotheses of Theorem \ref{T:di}. Moreover, in this case we can get better
estimates on the integers $j_r$ from Theorem \ref{T:di}.
\begin{Corollary}\label{C:sc1di} Let $f$ be an entire function of
order $\rho$ which satisfies (\ref{e:sc1}) for all $r$
sufficiently large. Then there is a constant $a>1$ such that
$Z_n(r)\le a(nm(ar)+n^2)$, for all $n\geq1$ and $r\geq1$.
\end{Corollary}
\begin{pf} Fix $r\geq1$, let $I_j=I(R_j,\rho/2,\beta,1,C)$ be the
covering system of admissible intervals from the proof of Theorem \ref{T:sc1},
and recall that $R_j\in((2k)^j,(2k)^{j+1})$. By Corollary
\ref{C:nz} we have $Z_n(r)\le 10nm(4R_{j_r})+An^2$ for all $n\geq1$, where
$A$ is a constant and $j_r$ is defined in (\ref{e:jr}) as the smallest integer
such that $R_{j_r}>rR_1(f)$ and $\Lambda m(4R_{j_r})\geq4CC_1m(4r)$.
\par We fix $j_0,j_1$ such that
$$(2k)^{j_0}\geq R_1(f),\;A_1^{j_0}\geq4CC_1/\Lambda,\;
r\in\left[(2k)^{j_1},(2k)^{j_1+1}\right).$$ Then
$R_{j_0+j_1+1}>rR_1(f)$ and
$$m(4R_{j_0+j_1+1})\geq A_1^{j_0}m\left(4(2k)^{j_1+1}\right)
>\frac{4CC_1}{\Lambda}\,m(4r).$$ Consequently, $j_r\leq j_0+j_1+1$,
$R_{j_r}\leq(2k)^{j_0+j_1+2}\leq(2k)^{j_0+2}r$,
and the corollary follows.\end{pf}
\par Given an entire function $f$, it is frequently known that $f$
verifies a growth condition $m(r)\le r^{\phi(r)}$, where
$\lim_{r\to\infty}\phi(r)=\rho$ and the function
$r^{\phi(r)-\rho}$ is slowly increasing. In the remainder of this
section, we denote by $r_n$ the unique solution of the equation
$r^{\phi(r)}=n$. Our next theorem shows that in this case there
are estimates $e_n=O(n^\tau\log n)$, provided that
$$m(r_{n_j})\geq ar_{n_j}^{\phi(r_{n_j})}=an_j$$ holds for
a ``slow growing" subsequence $n_j$.
\begin{Theorem}\label{T:pesif} In the above setting, assume there is an increasing
sequence of integers $n_j$ such that $n_{j+1}^\gamma\le bn_j$ and
$m(r_{n_j})\ge an_j$, where $0<\gamma\le1,\,b>0,\,0<a\le1$. Then
there exists a sequence of positive $\epsilon_n\rightarrow0$, such
that the estimates
\begin{eqnarray*}e_n&\leq&\frac{4C(3M)^{1/\gamma}}{\rho\gamma\Lambda^{1/\gamma}}\,
n^{1+1/\gamma}\log\frac{3Mn}{\Lambda}\,,\\
m_n(r)&\leq&\frac{3C(3M)^{1/\gamma}}{\Lambda^{1/\gamma}}\,n^{1+1/\gamma}\log
r,\;1\le r\le \frac12\,n^{1/\rho-\epsilon_n},\end{eqnarray*} hold
for all $n$ sufficiently large, where
$$C=\frac{2^{3\rho+4}}a\left(1+\frac1\rho\,\log_2(32/a)\right),\;\;
M=\frac{2^{(2\rho+3)\gamma}Cb}{a}\,.$$
\end{Theorem}
\begin{pf} We let $s_j=r_{n_j}$.
By Lemma \ref{L:si} with $k=(32/a)^{1/\rho}$ and $j$ sufficiently
large, there is $R_j\in(s_j/k,2s_j)$ such that
$$CS(R_j)\geq m(8s_j)\ge m(4R_j),$$ where
$$C=\frac{2^{3\rho+4}}a\left(1+\frac1\rho\,\log_2(32/a)\right).$$
We may assume that
$$(4s_j)^{\phi(4s_j)}\le
2\cdot 4^\rho s_j^{\phi(s_j)}=2^{2\rho+1}n_j.$$
Using this and (\ref{e:sm}) we get
$$S(R_j)\leq\frac{2m(2R_j)}{\log 2}\leq
4m(4s_j)\leq2^{2\rho+3}n_j.$$ Since $CS(R_j)\geq m(8s_j)\geq an_j$
it follows that
$$\frac{S^\gamma(R_{j+1})}{S(R_j)}\le
\frac{2^{(2\rho+3)\gamma}Cn_{j+1}^\gamma}{an_j}\le
\frac{2^{(2\rho+3)\gamma}Cb}{a}=M.$$ Moreover,
$$S(R_j)\ge\frac{a}{C}\,s_j^{\phi(s_j)}
\geq\frac{a}{C}\,(R_j/2)^{3\rho/4}\geq R_j^{\rho/2},$$ when $j$ is
sufficiently large, and if $\beta=\Lambda/(3M)$ then
$$\max\{\beta S^\gamma(R_j),\beta S^\gamma(R_{j+1})\}\leq\frac{\Lambda}{3}\,
S(R_j)\leq\frac{\Lambda}{2}\,S(R_j)-n_0-1.$$ So the intervals
$I_j=I(R_j,\rho/2,\beta,\gamma,C)$ form a covering system of
admissible intervals, starting with some $j$ sufficiently large.
The theorem now follows from Corollary \ref{C:me2}. If $n\in
I_{j_n}$ then the fact that $R_{j_n}\ge n^{1/\rho-\epsilon_n}/2$
can be proved exactly like the similar statement in Theorem
\ref{T:liminf}.
\end{pf}
\par This theorem has a corollary which allows us to estimate
$e_n$ and $m_n$ using the behavior of the Taylor coefficients of $f$.
\begin{Corollary}\label{C:tc} Suppose that for an entire function
$$f(z)=\sum_{n=0}^\infty c_nz^n$$ we have $m(r)\le r^{\phi(r)}$, where
$\lim_{r\to\infty}\phi(r)=\rho$ and the function
$r^{\phi(r)-\rho}$ is slowly increasing. Let $r_n$ be defined by
$r_n^{\phi(r_n)}=n$. If there is an increasing sequence of
integers $n_j$ such that $$n^\gamma_{j+1}\le bn_j\;,\;\;
\log|c_{n_j}|\ge an_j-n_j\log r_{n_j},$$ where
$0<\gamma\le1,\,b>0,\,0<a\leq1$, then the estimates on $e_n$ and
$m_n(r)$ from Theorem \ref{T:pesif} hold for all $n$ sufficiently
large.\end{Corollary}
\begin{pf} This follows from Theorem \ref{T:pesif}, since by the Cauchy
inequalities we have $m(r_{n_j})\ge\log|c_{n_j}|+
n_j\log r_{n_j}\ge an_j$.\end{pf}
\par As an example we take the entire function
$$f(z)=\sum_{j=1}^\infty(z/n_j)^{n_j},$$
where $n_1\geq2$, $n_{j+1}=n_j^{\tau-1}$ and $\tau>2$. This
function was studied in Section 6 of \cite{CP2}, where it was
shown that there are constants $C_1$ and $C_2$ such that $e_n\le
C_1n^\tau\log n$ and $m_n(r)\le C_2n^\tau\log r$ for $1\le r\le
n$. This was a result of quite elaborate estimates. Since $f$ is a
function of order 1 and type $1/e$, we have $m(r)\le 2r/e$ for $r$
large. Taking $r_{n_j}=en_j/2$ and $a=\log(e/2)$ we get
$$\log c_{n_j}=-n_j\log n_j=an_j-n_j\log r_{n_j}.$$
So Corollary \ref{C:tc} applies with $\gamma=1/(\tau-1)$ and
$b=1$.

\section{The functions $\zeta$ and $\xi$}\label{S:rf}
The Riemann $\zeta$-function is holomorphic in ${\Bbb C}$ except at
$z=1$, where it has a simple pole (see e.g. \cite[Theorem 2.1]{T}). The
function $\xi$ is defined by
$$\xi(z)=\frac{z(z-1)}{2}\,\pi^{-z/2}\Gamma\left(\frac z2\right)\zeta(z),$$
where $\Gamma$ is the Euler Gamma function (see \cite[(2.1.12)]{T}).
Then $\xi$ is an entire function of order 1 \cite[Theorem 2.12]{T}.
\par For the convenience of the reader, we include the proof of
the following proposition.
\begin{Proposition}\label{P:zeta} There exist positive constants $c_1<c_2$, $d_1<d_2$
such that for all $r\geq2$ we have
\begin{eqnarray*}&&c_1r\log r\leq m(r,\zeta)\leq c_2r\log r,\\
&&d_1r\log r\leq m(r,\xi)\leq d_2r\log r.
\end{eqnarray*}\end{Proposition}
\begin{pf} If $x={\mathbf {Re\,}} z>0$, then
$$\Gamma(z)=\int_0^{\infty}t^{z-1}e^{-t}\,dt\;,\;\;
|\Gamma(z)|\le\Gamma(x).$$ We let $\mu(r)$ be the supremum of
$\log|\Gamma(z)|$ when $|z|=r$ and $x\geq1/2$. Then by Stirling's
formula $C_1r\log r\le\mu(r)\le C_2r\log r$ for $r\geq2$.
\par For $x>0$ one has (see \cite[(2.1.4)]{T})
$$\zeta(z)=z\int_1^\infty\frac{[t]-t+\frac12}{t^{z+1}}\,dt+\frac1{z-1}+
\frac12\,.$$ Hence for $x\ge1/2$ and $|z-1|>2$ we have
$$|\zeta(z)|\le\frac{|z|}2\int_1^\infty\frac{1}{t^{x+1}}\,dt+1\le
|z|+1,$$ while $\zeta(x)=\sum_{k=1}^\infty k^{-x}>1$ for $x>1$. To
estimate $\zeta(z)$ for ${\mathbf {Re\,}} z\le1/2$ we use the
functional equation (see \cite[Theorem 2.1]{T})
$$\zeta(z)=2^z\pi^{z-1}\sin\frac{\pi z}2\,\Gamma(1-z)\zeta(1-z).$$
We conclude that $m(r,\zeta)\le c_2r\log r$. But for odd
integers $n>0$ we have $|\zeta(-n)|\ge2^{-n}\pi^{-n-1}n!$ and
therefore $m(n,\zeta)\ge c'_1n\log n$. This implies $m(r,\zeta)\ge
c_1r\log r$, $c_1>0$.
\par Using the definition of $\xi(z)$ we have for $|z|=r$ with ${\mathbf {Re\,}} z\ge1/2$
$$|\xi(z)|\le\frac{r(r+1)^2}{2}\,\Gamma\left(\frac{r}2\right),\;
\xi(r)\ge\frac{r(r-1)}2\,\pi^{-r/2}\Gamma\left(\frac r2\right).$$
Since $\xi(z)=\xi(1-z)$ (see \cite[(2.1.13)]{T}) we obtain for
$|z|=r$ with ${\mathbf {Re\,}} z\le1/2$
$$|\xi(z)|\le\frac{(r+1)(r+2)^2}{2}\,\Gamma\left(\frac{r+1}2\right).$$
So we see that $d_1r\log r\le m(r,\xi)\le d_2r\log r$.\end{pf}
\par By Proposition \ref{P:zeta} the function $\xi$ verifies condition
(\ref{e:sc1}) if $k$ is chosen sufficiently large, so Theorem
\ref{T:sc1} and Corollary \ref{C:sc1di} hold for $\xi$.
Since $\zeta$ is meromorphic, the quantity $m_n(r)$ is not well
defined for $\zeta$. We have the following:
\begin{Theorem}\label{T:zeta} There exists a constant $C>0$ such that for every
integer $n\geq1$ and every $r\ge1$ we have
$$\frac{n^2+3n}{2}\leq Z_n(r,\zeta)\leq C(nr\log r+n^2).$$ \end{Theorem}
\begin{pf} Since $\zeta$ is holomorphic near 0, we can find,
by a simple dimension argument, a non-trivial
polynomial $P(z,w)$ of degree at most $n$ such that $P(z,\zeta(z))$ has a
zero of order at least $(n^2+3n)/2$ at 0 (see the proof of Theorem 2.5
in \cite{CP2}). This implies the lower estimate on $Z_n$.
\par The function $\tilde\zeta(z)=(z-1)\zeta(z)$ is entire.
Proposition \ref{P:zeta} implies that
$$c_1r\log r\leq m(r,\tilde\zeta)\leq c'_2r\log r.$$ By
Corollary \ref{C:sc1di} it follows that there exists a constant
$C>0$ such that
$$Z_n(r,\tilde\zeta)\leq  C(nr\log r+n^2),$$ for all
$n\geq1$ and $r\ge1$. Note that if $P(z,w)$ is a polynomial of
degree at most $n$, then there exists a polynomial $Q(z,w)$ of
degree at most $2n$ such that
$(z-1)^nP(z,\zeta(z))=Q(z,\tilde\zeta(z))$. Hence
$Z_n(r,\zeta)\leq Z_{2n}(r,\tilde\zeta)$, and the proof is
complete.
\end{pf}

\section{Extremal functions}\label{S:ext}
\par If $K\subset{\Bbb C}^2$ is a compact set, the extremal function
$V_K$ of $K$ (also called the pluricomplex Green function of $K$
with pole at infinity) is defined by
$$V_K(z,w)=\sup\,\frac{1}{\deg P}\log|P(z,w)|,$$
where the supremum is taken over all polynomials $P$ such that
$\|P\|_K\leq1$. Then either $V_K$ is finite at every point, or
$V_K\equiv\infty$, and the latter occurs if and only if $K$ is
pluripolar (see e.g. \cite[Ch. 5]{K}).
\par Let $f$ be an entire transcendental function and let $$K=\{(z,f(z)):\,|z|\leq1\}.$$
Then $K$ is pluripolar and $V_K\equiv\infty$. Using our estimates
on $m_n(r)$, it is still possible to define a meaningful extremal
function of $K$ along the graph of $f$. This relates to
Sadullaev's result on the existence of extremal functions for
non-pluripolar subsets of algebraic varieties \cite{Sa}.
\par We assume in this section that $f$ is an entire transcendental function
which verifies $$m_n(r)\leq C_fn^2\log r,\;1\leq r\leq
r_n,\;n\geq1,$$ where $r_n$ is a sequence increasing to infinity
and $C_f$ is a constant depending on $f$. Classes of such
functions are constructed in Section \ref{S:scof}. Let us define
$$W_n(z)=\sup\log|P(z,f(z))|\,,$$
where the supremum is taken over all polynomials $P$ of degree at
most $n$ which verify $|P(z,f(z))|\leq1$ on $\Delta$. The
functions $W_n$ are non-negative, continuous and subharmonic on
${\Bbb C}$, and $W_n\equiv0$ on $\Delta$. \par Next we define
$$W(z)=\limsup_{n\rightarrow\infty}\frac{1}{n^2}\,W_n(z),$$
and we let $W^\star$ denote the upper semicontinuous
regularization of $W$. We have the following:
\begin{Proposition}\label{P:extremal}The function $W^\star$ is non-negative
subharmonic on ${\Bbb C}$, $W^\star\equiv0$ on $\Delta$, and for
all $r\geq1$
$$\frac{1}{2}\,\log r\leq\max\{W^\star(z):\,|z|=r\}\leq C_f\log r.$$
If $f(z)=e^z$ then $W^\star(z)=\frac{1}{2}\,\log^+|z|$ for all
$z\in{\Bbb C}$.\end{Proposition}
\begin{pf} By the proofs of Theorem 2.5 and Corollary 2.6 of \cite{CP2}
there exists, for each $n\geq1$,
a non-trivial polynomial $P_n(z,w)$ of degree $n$, such that the
function $F_n(z)=P_n(z,f(z))$ verifies $M(1,F_n)=1$ and
$$\frac{n^2+3n}{2}\,\log r\leq\log M(r,F_n)\leq m_n(r),$$ for all
$r\geq1$. Note that, in particular, this implies $C_f\geq1/2$ for
any $f$.
\par Our assumption on the growth of $m_n(r)$ implies that the family of subharmonic functions
$W_n/n^2$ is locally upper bounded, hence $W^\star$ is a
non-negative subharmonic function on ${\Bbb C}$ which verifies
$$W^\star(z)\leq C_f\log^+|z|\,.$$ Suppose that for some $r>1$ we
have
$$\max\{W^\star(z):\,|z|=r\}<\frac{1}{2}\,\log r.$$ The Hartogs Lemma implies that for
$n$ large and for all $z$ with $|z|=r$
$$\log|F_n(z)|\leq W_n(z)<\frac{n^2}{2}\,\log r.$$ This contradicts the above lower estimate
on $\log M(r,F_n)$. \par In the case of the function $f(z)=e^z$ it
was proved in \cite{CP1} that
$$m_n(r)\leq\left(\frac{n^2}{2}+o(n^2)\right)\log r,\;1\leq r\leq
n.$$ The preceding argument shows that now
$W^\star(z)\leq\frac{1}{2}\,\log^+|z|$. We conclude that the
equality must hold, by applying the maximum principle on ${\Bbb
C}\setminus\Delta$ to the subharmonic function
$W^\star(z)-\frac{1}{2}\,\log|z|\leq0$.\end{pf}
\section{Estimates for algebraic measures}\label{S:eam}
\par Throughout Sections \ref{S:eam} and \ref{S:agtf},
$K$ is an algebraic extension of degree $\sigma$ of the field
${\Bbb Q}$ of rational numbers and $f$ is, unless otherwise
specified, an entire transcendental function of finite positive
order $\rho$. Without loss of generality we may assume that
$M(r,f)\geq r$ for $r\geq1$.
\par For an algebraic number $\zeta$, we define its norm
$\|\zeta\|$ as the maximum of the absolute values of its
conjugates. Then $\|\zeta_1\zeta_2\|\le \|\zeta_1\|\|\zeta_2\|$
and $\|\zeta_1+\zeta_2\|\le \|\zeta_1\|+\|\zeta_2\|$ (see \cite[p.
62]{M}).
\par If $P(\zeta_1,\dots,\zeta_n)$ is a polynomial with algebraic
coefficients, then its {\it height} $h(P)$ is defined as the
maximum of the norms of its coefficients.
\par If $\omega_1,\dots,\omega_\sigma$ is a basis for the ring $I_K$ of algebraic
integers in $K$, then any  $\zeta\in I_K$ can be written as
\begin{equation}\label{e:am}
\zeta=p_1\omega_1+\dots+p_\sigma\omega_\sigma,\end{equation} where
$p_1,\dots,p_\sigma$ are rational integers. If
$|||\zeta|||=\max\{|p_1|,\dots,|p_\sigma|\}$, then (see \cite[p.
62]{M}) there are constants $\gamma_1$ and $\gamma _2$ depending
only on $K$ such that
$\gamma_1|||\zeta|||\le\|\zeta\|\le\gamma_2|||\zeta|||$.
\par Given a natural number $d$ we denote by $I_K(d)$ the set of
numbers $z\in K$ such that $dz\in I_K$, and by $I_K(d,A)$ the set
of $z\in I_K(d)$ with $\|z\|\le A$. Let $N_K(d,A,r)$ be the number
of points in $I_K(d,A)\cap\Delta_r$. \begin{Lemma}\label{L:np}
There exist constants $c'$ and $c''$ depending only on $K$ with
the following properties: If $K$ is real and $r\leq A$ then
$$c'd^\sigma A^{\sigma-1}r\le N_K(d,A,r)\le c''d^\sigma
A^{\sigma-1}r.$$ If $K$ is complex and $r\leq A$ then $$c'd^\sigma
A^{\sigma-2}r^2\le N_K(d,A,r)\le c''d^\sigma A^{\sigma-2}r^2.$$
Moreover, in both cases, if $r>A$ then
$$c'(dA)^\sigma\le N_K(d,A,r)\le c''(dA)^\sigma.$$ \end{Lemma}
\begin{pf}  We will consider only the complex case, so $\sigma\ge 2$.
The real case can be considered in a similar manner.
\par We assume at first that $d=1$. If $N$ is the number
of points in $I_K(1,A)$, then by (\ref{e:am}) $c_1A^\sigma\le N\le
c_2A^\sigma$, where the constants $c_1$ and $c_2$ depend only on
$K$. Moreover $I_K(1,A)\subset\Delta_A$, so the lemma is proved
for $r>A$.
\par Suppose that $r\le A$. There exists an absolute constant
$c_3$ such that $\Delta_A$ can be covered by $c_3A^2r^{-2}$ disks
of radius $r/2$. Therefore there is a disk $D$ of radius $r/2$
containing at least $c_4A^{\sigma-2}r^2$ points of $I_K(1,A/2)$,
where $c_4$ depends only on $K$. Let $z_0$ be one of this points.
If $z\in I_K(1,A/2)\cap D$, then $z-z_0\in
I_K(1,A)\cap\Delta_{r}$. Consequently, $N_K(1,A,r)\ge
c'A^{\sigma-2}r^2$, where $c'$ depends only on $K$.
\par Now let $N=N_K(1,A,r)$ and let $\omega_1,\dots,\omega_\sigma$ be
a basis for $I_K$ over ${\Bbb Z}$. Since $I_K\not\subset{\Bbb R}$,
we may assume that $\omega_1/\omega_2\not\in{\Bbb R}$. Then it is
easy to see that there is a constant $c_5$ depending only on $K$
such that the disk $\Delta_A$ contains at least $c_5A^2$ points
$z$ of the form $z=p_1\omega_1+p_2\omega_2$, where
$p_1,p_2\in{\Bbb Z}$ and $|p_1|,|p_2|\le A/\gamma_2$, so
$\|z\|\leq A$. Moreover, there is a constant $c_6$ depending only
on $K$  and at least $c_6A^2r^{-2}$ disjoint disks of radius $r$
centered at these points. Hence each of these disks contains at
least $N$ points from $I_K(1,2A)$. It follows that
$c_6NA^2r^{-2}\le c_22^\sigma A^\sigma$, so $N\le
c''A^{\sigma-2}r^2$.
\par If $d>1$ then we note that $z\in I_K(1,dA,dr)$ if and only if
$z/d\in I_K(d,A,r)$, hence $N_K(d,A,r)=N_K(1,dA,dr)$.
\end{pf}
\par We say that a function $f$ takes values at $z$ in $I_K(d)$ with
multiplicity $m$ if the numbers $z, f(z),\dots,f^{(m-1)}(z)$
belong to $I_K(d)$. In this case we define $\|f(z)\|_m$ as the
maximum of $\|z\|, \|f(z)\|,\dots,\|f^{(m-1)}(z)\|$.
\par In this setting, we have the following lemma (see Ch. 1, \S 2
and Ch. 2, \S 2 in \cite{G}).
\begin{Lemma}\label{L:gt} Let $f$ be a holomorphic function in
a neighborhood of $z_0$, which takes values at $z_0$ in $I_K(d)$
with multiplicity $m$, and such that $\|f(z_0)\|_m\le A$, $A\ge1$.
If $P(z,w)$ is a polynomial of degree $n$ with coefficients in
$I_K$ and of height $h$, and if $F(z)=P(z,f(z))$, then for $k\leq
m-1$ we have $F^{(k)}(z_0)\in I_K(d^n)$ and
$$\|g^{(k)}_{ij}(z_0)\|\le A^{i+j}(i+j)^k,$$
where $g_{ij}(z)=z^if^j(z)$, $i+j\leq n$. Moreover, if
$F^{(k)}(z_0)\neq0$ then
$$|d^nF^{(k)}(z_0)|\ge
\left(hd^nA^{n}(n+1)^{k+2}\right)^{-\sigma+1}.$$ \end{Lemma}
\begin{pf} If
$$f(z)=\sum_{j=0}^\infty a_{j}(z-z_0)^j,$$
then the $k$-th derivative $a_{jk}$  of $f^j(z)$ at $z_0$ is
$$a_{jk}=k!\sum_{i_1+\dots+i_j=k}a_{i_1}\cdots a_{i_j}=
k!\sum_{i_1+\dots+i_j=k}\frac{f^{(i_1)}(z_0)\cdots
f^{(i_j)}(z_0)}{i_1!\cdots i_j!}.$$ Hence $d^ja_{jk}\in I_K$.
Since
$$j^k=\sum_{i_1+\dots+i_j=k}\frac{k!}{i_1!\cdots i_j!}\,,$$
we see that $\|a_{jk}\|\le j^kA^j$.
\par If $g_{ij}(z)=z^if^j(z)$, then
$$g^{(k)}_{ij}(z_0)=\sum_{p=0}^{\min\{i,k\}}{k\choose p}
\frac{i!}{(i-p)!}z_0^{i-p}a_{j,k-p}\,.$$ Thus
$d^{i+j}g^{(k)}_{ij}(z_0)\in I_K$. Moreover,
\begin{eqnarray*}\|g^{(k)}_{ij}(z_0)\|&\leq&
\sum_{p=0}^{\min\{i,k\}}{k\choose p}\frac{i!}{(i-p)!}j^{k-p}A^{i+j-p}\\
&\leq&A^{i+j}\sum_{p=0}^k{k\choose
p}i^pj^{k-p}=A^{i+j}(i+j)^k.\end{eqnarray*} Hence
$d^nF^{(k)}(z_0)\in I_K$ and
$$\|F^{(k)}(z_0)\|\le\frac{(n+1)(n+2)}{2}\,hA^{n}n^k\leq hA^n(n+1)^{k+2}.$$
\par Since $c_1=d^nF^{(k)}(z_0)\in I_K$, the number $\nu$ of its conjugates
$c_2,\dots,c_\nu$ does not exceed $\sigma$ and
$$\left|\prod_{i=1}^\nu c_i\right|\ge1$$
when $c_1\ne0$.  Note that $|c_i|\le\|c_i\|=\|c_1\|$.
Consequently,
$$|c_1|\ge \left(hd^nA^{n}(n+1)^{k+2}\right)^{-\sigma+1}.$$
\end{pf}
\par The following result is a
consequence of C. L. Siegel's lemma adapted for our purposes.
\begin{Lemma}\label{L:vsl} Suppose that there are $l$ points $z_1,\dots,z_l$
in $\Delta_r$, $r\ge1$, such that, for $1\le q\le l$, $f$ takes
values at $z_q$ in $I_K(d_q)$ with multiplicity $m_q$ and
$\|f(z_q)\|_{m_q}\le A$, $A\ge1$. If $\nu=\sum_{q=1}^lm_q<N$,
where $N=(n+1)(n+2)/2$, and $m_q\le m$, $d_q\leq d$, then there is
a polynomial $P(z,w)$ of degree $n$ with coefficients $c_{ij}\in
I_K$ and of height
$$h(P)\leq H_n=C_1\left(C_2d^nA^n(n+1)^{m+1}\right)^{\nu/(N-\nu)},$$
where $C_1,C_2$ are constants depending only on $K$, with the
following properties: The function
$$F(z)=P(z,f(z))=\sum_{i+j=0}^nc_{ij}z^if^j(z)\not\equiv0$$
and for $t\ge2r$
$$\|F\|_{\Delta_{2r}}\le
(n+1)^2H_n\left(\frac{4r}{t}\right)^\mu M^n(t,f),$$ where
$\mu\ge\nu$ is the number of zeros of $F$ in
$\Delta_r$.\end{Lemma}
\begin{pf} By Lemma \ref{L:gt} $d_q^ng^{(k)}_{ij}(z_q)$, where
$g_{ij}(z)=z^if^j(z)$, is an algebraic integer and
$$\|d_q^ng^{(k)}_{ij}(z_q)\|\leq d_q^nA^{n}n^k\leq d^nA^nn^{m-1}.$$
\par Let us consider the system of $\nu$ equations
$$\sum_{i+j=0}^nc_{ij}d_q^ng_{ij}^{(k)}(z_q)=0,\qquad 1\le q\le l,
\qquad 0\le k\le m_q-1,$$ with $N$ unknowns $c_{ij}$. By \cite[p.
63]{M} there are constants $C_1$ and $C_2$ depending only on $K$
such that this system has a non-trivial solution in $I_K$ with
$$\|c_{ij}\|\leq
C_1\left(C_2Nd^nA^nn^{m-1}\right)^{\nu/(N-\nu)}\leq H_n.$$
\par Since $\|P\|_{\Delta^2}\le (n+1)^2H_n$, by the Bernstein--Walsh
inequality
$$|P(z,w)|\le (n+1)^2H_n\exp(n\max\{\log^+|z|,\log^+|w|\}),$$
so $\|F\|_{\Delta_t}\le (n+1)^2 H_nM^n(t,f)$.
\par The function $F$ has  $\mu\ge\nu$ zeros in $\Delta_r$, so by
Lemma \ref{L:bp}
$$\|F\|_{\Delta_{2r}}\le \left(\frac{4r}{t}\right)^\mu\|F\|_{\Delta_t}
\leq(n+1)^2H_n\left(\frac{4r}{t}\right)^\mu M^n(t,f).$$\end{pf}
\par Suppose that for a set $E\subset{\Bbb C}$ and for some integer
$m\geq1$ we have $z,f(z),\dots,f^{(m-1)}(z)\in K$ for all $z\in
E$. Then for $z\in E$ we let $d_z$ be the smallest natural number
such that $f$ takes values at $z$ in $I_K(d_z)$ with multiplicity
$m$. We set
\begin{eqnarray*}
&\|f\|_{E,m}=\max\left\{1,\sup_{z\in
E}\|f(z)\|_m\right\},\;d(E,m)=\sup_{z\in E}d_z,\\
&{\mathcal A}_K(E,m)=d(E,m)\|f\|_{E,m}.
\end{eqnarray*}
The number ${\mathcal A}_K(E,m)$ will be called the {\it algebraic
measure of order $m$} of the function $f$ on $E$. If for some
$z\in E$ we have $z\not\in K$ or $f^{(j)}(z)\not\in K$ for some
$j<m$, then we set ${\mathcal A}_K(E,m)=\infty$. Note also that if
a set $E$ is infinite then ${\mathcal A}_K(E,m)=\infty$ for every
$m\geq1$.
\par Throughout the rest of this section and in Section \ref{S:agtf}
we will assume that $m(r)\le r^{\phi(r)}$, where
$\lim_{r\to\infty}\phi(r)=\rho$ and the function
$r^{\phi(r)-\rho}$ is slowly increasing. We denote by $r_n$ the
unique solution of the equation $r^{\phi(r)}=n$. The following
result is the main tool in the forthcoming estimates of the
algebraic measure.
\begin{Theorem}\label{T:geap} There exists a constant $C_K$ depending only on $K$
with the following property: If $n\geq1$, $1\le r\le r_n/4$ and
$E\subset\Delta_r$, then there are integers $k\geq0$ and $\mu$
such that
\begin{eqnarray*}&&(k+1)|E|>n^2/4,\;\max\{n^2/4,k|E|\}\le\mu\leq Z_n(r),\\
&&C_K{\mathcal A}_K^{2\sigma}(E,k+1)\geq\left(\frac
r{k(n+1)^{2\sigma-1}}\right)^{k/n} \exp\left(\frac
\mu{n}\log\frac{r_n}{4e^4r}\right).\end{eqnarray*}
\end{Theorem}
In the above statement we let $k^k=1$ if $k=0$.
\begin{pf} We may assume that $E$ is finite. Let $E=\{z_1,\dots,z_l\}$
and $\nu=[n^2/4]+1$. Note that $\nu/(N-\nu)\leq1$ and by Theorem
2.5 and Corollary 2.6 in \cite{CP2} we have $\nu\leq(n^2+3n)/2\leq
Z_n(r)$.
\par Let $m=[\nu/l]$. If ${\mathcal A}_K(E,m+1)=\infty$, then
we take $k=m$ and $\mu=\nu$ and the proof is finished. Otherwise,
we let $A_1=\|f\|_{E,m+1}$ and $d_1=d(E,m+1)$. We have $\nu=ml+p$,
$0\le p\le l-1$.  Applying Lemma \ref{L:vsl} with the above
points, with $m_q=m+1$ when $1\le q\le p$ and $m_q=m$ when $p+1\le
q\le l$, and with this value of $\nu$, we construct a non-trivial
polynomial $P(z,w)$ of degree $n$ with coefficients in $I_K$ and
with height
$$h(P)\le h=C_1C_2d_1^{n}A_1^n(n+1)^{m+2}$$ such that
$$\|F\|_{\Delta_{2r}}\le
(n+1)^2h\left(\frac{4r}{t}\right)^\mu M^n(t,f),$$ where $\mu\ge
n^2/4$ is the number of zeros of $F$ in $\Delta_r$ and $t\ge2r$.
\par There exist $q$ and $k$, $1\le q\le l$,
$0\le k\le\mu/l$, such that $F^{(k)}(z_q)\ne0$. We may assume that
${\mathcal A}_K(E,k+1)<\infty$. Clearly, $k\ge m$ and therefore
$kl>\nu-l>n^2/4-l$. Moreover $A=\|f\|_{E,k+1}\ge A_1$,
$d=d(E,k+1)\ge d_1$, so
$$h\le C_1C_2d^{n}A^n(n+1)^{k+2}.$$
By Lemma \ref{L:gt}
\begin{eqnarray*}
|d^nF^{(k)}(z_q)|&\ge&\left(hd^nA^n(n+1)^{k+2}\right)^{-\sigma+1}\\
&\geq&\left(C_1C_2d^{2n}A^{2n}(n+1)^{2k+4}\right)^{-\sigma+1}.
\end{eqnarray*}
\par By Cauchy's inequalities
\begin{eqnarray*}|d^nF^{(k)}(z_q)|&\le
&d^nk!\frac{\|F\|_{\Delta_{2r}}}{r^{k}} \leq d^n
(n+1)^2h\left(\frac kr\right)^k
\left(\frac{4r}{t}\right)^{\mu}M^n(t,f)\\
&\leq& C_1C_2d^{2n}A^n(n+1)^{k+4}\left(\frac
kr\right)^k\left(\frac{4r}{t}\right)^{\mu}M^n(t,f).
\end{eqnarray*}
We obtain
\begin{eqnarray*}
&&\left(C_1C_2d^{2n}A^{2n}(n+1)^{2k+4}\right)^{-\sigma+1}
\leq\\
&&C_1C_2d^{2n}A^n(n+1)^{k+4}\left(\frac
kr\right)^k\left(\frac{4r}{t}\right)^{\mu}M^n(t,f).
\end{eqnarray*} Since $A\geq1$ this implies
$$(C_1C_2)^\sigma(d^{2n}A^{2n})^\sigma\geq
\left(\frac rk\right)^k
\left(\frac{t}{4r}\right)^{\mu}M^{-n}(t,f)(n+1)^{-k(2\sigma-1)-4\sigma}.$$
Let $C_K=(16C_1C_2)^\sigma$. Taking the $n$th root and using the
inequality $(n+1)^{-1/n}\geq1/2$, we get
$$C_K{\mathcal A}^{2\sigma}(E,k+1)\geq
\left(\frac r{k(n+1)^{2\sigma-1}}\right)^{k/n}
\left(\frac{t}{4r}\right)^{\mu/n}M^{-1}(t,f).$$
\par Let $t=r_n$. Since $\mu\ge n^2/4$, $4r\leq r_n$ and $M(r_n,f)\leq e^n$,
we get $$\left(\frac{r_n}{4r}\right)^{\mu/n}M^{-1}(r_n,f)\ge
\left(\frac{r_n}{4r}\right)^{\mu/n}e^{-4\mu/n}=\exp\left(\frac
{\mu}{n}\log\frac{r_n}{4e^4r}\right),$$ so
$$C_K{\mathcal A}_K^{2\sigma}(E,k+1)\geq\left(\frac r{k(n+1)^{2\sigma-1}}\right)^{k/n}
\exp\left(\frac{\mu}{n}\log\frac{r_n}{4e^4r}\right).$$
\end{pf}
\section{Algebraic growth of transcendental functions}\label{S:agtf}
\par Let $f$ be an entire function of finite order $\rho$ and $K$
be an algebraic number field of degree $\sigma=[K:{\Bbb Q}]$. As
in Section \ref{S:eam} we assume, without loss of generality, that
$m(r)\leq r^{\phi(r)}$, where $\phi(r)\rightarrow\rho$ and
$r^{\phi(r)-\rho}$ is a slowly increasing function. Recall the
definition of the sequence $\{r_n\}$ by the equations
$r_n^{\phi(r_n)}=n$.
\par Given a transcendental function $f$ we define the {\it algebraic
growth characteristic} of $f$ on $K$ by
$${\mathbf a}_K(s,r,m)=\inf\{\log {\mathcal A}_K(E,m):\,E\subset\Delta_r, |E|\ge
s\}.$$ Due to our knowledge of the behavior of $Z_n(r)$ we are now
able to get estimates for ${\mathbf a}_K(s,r,m)$. The first series
of results applies to general transcendental functions. We recall
that when $\sigma>2$, then for every $\epsilon>0$ there is an
entire function $f$ of order smaller than $\epsilon$ such that
$f(K)\subset I_K$ (see \cite[Satz 1]{GS}) . Moreover, one can find
such a function so that $f^{(m)}(K)\subset K$ for all $m$ (see
\cite[Satz 2]{GS}) .
\par Our first theorem shows that when $m$ and $r$ are fixed, the
algebraic growth characteristic exceeds $s^{1/2}\log s$, at least
for a subsequence of integers $s$.
\begin{Theorem}\label{T:id}If $f$ has finite order $\rho>0$ then for all $m,r\geq1$
$$\limsup_{s\to\infty}\frac{{\mathbf a}_K(s,r,m)}{s^{1/2}\log s}>
\frac{2^{-3\rho/2-10}}{\sigma}\left(\frac{\Lambda
m}{\rho(\rho+5)}\right)^{1/2}.$$
\end{Theorem}
\begin{pf} By Corollary \ref{C:gnz} there is a fundamental sequence of
integers $\{n_j\}$ for $f$ with the following property: For every
$r\ge1$ there is an integer $j_r$ such that $Z_{n_{j}}(r)\le a
n_{j}^2$ for $j\geq j_r$, where
$$a=\frac{2^{3\rho+11}(\rho+5)}{\Lambda \rho}\,.$$ We may assume that
$4r\le r_{n_j}$ when $j\ge j_r$. Let $s_j=an_j^2/m+1$ and $E$ be a
subset of $\Delta_r$ with $|E|\ge s_j$. If $k$ is the integer from
Theorem \ref{T:geap} corresponding to $n=n_j,r,E$, then
$m|E|>Z_{n_j}(r)\geq k|E|$, so $m\geq k+1$. It follows from
Theorem \ref{T:geap} that
$$C_K{\mathcal A}_K^{2\sigma}(E,m)\geq
\left({m(n_j+1)^{2\sigma-1}}\right)^{-m/n_j}
\exp\left(\frac {n_j}4\log\frac{r_{n_j}}{4e^4r}\right).$$
\par For all $j$ sufficiently large (depending on $r,m$) we have
$$\left(m(n_j+1)^{2\sigma-1}\right)^{-m/n_j}\ge1/2.$$ Since
$r_{n_j}^{\phi(r_{n_j})}=n_j$ and $\phi(r_{n_j})\rightarrow\rho$,
we conclude that there is a sequence of positive $\epsilon_j\to0$
such that
$$\frac {n_j}4\log\frac{r_{n_j}}{4e^4r}\geq
\frac{(1-\epsilon_j)n_j}{4\rho}\log n_j.$$ It follows that
$$2\sigma\,{\mathbf a}_K(s_j,r,m)+\log(2C_K)\geq
\frac{(1-\epsilon_j)n_j}{4\rho}\log n_j.$$ Since
$n_j=(m(s_j-1)/a)^{1/2}$ we see that
$$\limsup_{s\to\infty}\frac{{\mathbf a}_K(s,r,m)}{s^{1/2}\log s}\ge
\frac{(m/a)^{1/2}}{16\sigma\rho}>\frac{2^{-3\rho/2-10}}{\sigma}
\left(\frac{\Lambda m}{\rho(\rho+5)}\right)^{1/2}.$$
\end{pf}

\vspace{2mm}\noindent {\bf Remark:} It is interesting to note that
the value of $\limsup$ in the above theorem is achieved on a
sequence $\{s_j\}$ depending only on $f$ and $m$,
$s_j=an_j^2/m+1$.

\vspace{2mm}
\par As mentioned above, there are functions whose derivatives of all
orders map $K$ into $K$. So it is interesting to estimate the
algebraic growth of such functions on the sets $I_K(d,A)$. The
number $dA$ can be viewed as the algebraic measure of the set
$I_K(d,A)$, while ${\mathcal A}_K(I_K(d,A),m)$ is the algebraic
measure of the set of values of $f$ and its derivatives on
$I_K(d,A)$. We introduce
$$\eta_K(\lambda,r,m)=
\inf\{\log{\mathcal A}_K(I_K(d,A)\cap\Delta_r,m):\,dA\ge\lambda\}.$$
The following result describes the growth of
$\eta_K(\lambda,r,m)$.
\begin{Corollary}\label{C:iad} If $\sigma\ge3$, then there is a constant $c'$
depending only on $K$ such that for $r,m\geq1$
$$\limsup_{\lambda\to\infty}
\frac{\eta(\lambda,r,m)}{\lambda^{\sigma/2-1}\log
\lambda}\ge\frac{\sigma-2}{\sigma}\,
2^{-3\rho/2-10}\left(\frac{c'\Lambda
m}{\rho(\rho+5)}\right)^{1/2}.$$
\end{Corollary}
\begin{pf} By Lemma \ref{L:np}, $|I_K(d,A)\cap\Delta_r|\ge c'(dA)^{\sigma-2}$.
By Theorem \ref{T:id}, let $s_j$ be a sequence such that
$$\limsup_{j\to\infty}\frac{{\mathbf a}_K(s_j,r,m)}{s_j^{1/2}\log s_j}\ge
\frac{2^{-3\rho/2-10}}{\sigma}\left(\frac{\Lambda
m}{\rho(\rho+5)}\right)^{1/2}.$$ We define $\lambda_j$ by
$s_j=c'\lambda_j^{\sigma-2}$. If $dA\ge\lambda_j$ then
$\eta(\lambda_j,r,m)\ge{\mathbf a}_K(s_j,r,m)$, so
$$\limsup_{j\to\infty}
\frac{\eta(\lambda_j,r,m)}{\lambda_j^{\sigma/2-1}\log
\lambda_j}\geq\sqrt{c'}(\sigma-2)
\limsup_{j\to\infty}\frac{{\mathbf a}_K(s_j,r,m)}{s_j^{1/2}\log
s_j}\,,$$ and the conclusion follows.\end{pf}
\par Let $I_K(A)=I_K(1,A)$ be the set of algebraic integers $z\in I_K$
of norm $\|z\|\leq A$. Clearly, $I_K(A)\subset\Delta_A$. In our
next theorem we estimate the number of points $z\in I_K(A)$ which
are mapped to points of $I_K$ of smallest possible norm $A'$.
Since $|z|\leq A$ and $|f(z)|\leq\|f(z)\|$, it is natural to
expect, due to the growth of $f$, that $A'\geq\exp(A^{\phi(A)})$.
We will prove that if $\rho<\sigma/2$ then the proportion of
points of $I_K(A_j)$ which are mapped by $f$ into $I_K(\exp
A_j^{\phi(A_j)})$ tends to 0, for a certain sequence
$A_j\rightarrow\infty$.
\par To this end, we need the following version of Theorem
\ref{T:geap}, which provides upper bounds for $|E|$ if the
algebraic measure of order 1 of $f$ on $E$ is bounded above by
certain quantities. \begin{Proposition}\label{P:geap}There exists
a constant $C_K$ depending only on $K$ such that if $n\geq1$,
$1\leq r\leq r_n/4$, $E\subset\Delta_r$ and
$$C_K{\mathcal A}_K^{2\sigma}(E,1)<\exp\left(\frac{n}{4}\log\frac{r_n}{4e^4r}\right),$$
then $|E|\leq Z_n(r)$.\end{Proposition}
\begin{pf} If $|E|>Z_n(r)$ and $k$ is the integer from Theorem \ref{T:geap},
then $k|E|\leq Z_n(r)$ implies $k=0$. Since $\mu\geq n^2/4$, we
reach a contradiction with the conclusion of Theorem
\ref{T:geap}.\end{pf}
\par We have the following theorem:
\begin{Theorem}\label{T:proportion} If $f$ is an entire function of order $0<\rho<\sigma/2$ then
$$\liminf_{A\rightarrow\infty}
\frac{\left|I_K(A)\cap f^{-1}(I_K(\exp
A^{\phi(A)}))\right|}{|I_K(A)|}=0.$$
\end{Theorem}
\begin{pf} By Theorem \ref{T:liminf} and Corollary 2.6 in \cite{CP2}, there
exists a fundamental sequence $\{n_j\}$ for $f$ and positive
numbers $\epsilon_j\rightarrow0$ such that
$$Z_{n_j}(r)\leq an_j^2\log3r,\;1\leq r\leq n_j^{1/\rho-\epsilon_j}/6,$$
where $a=2^{3\rho+10}(\rho+5)/(\Lambda\rho)$.
\par Let $A_j=n_j^{1/((1+\epsilon)\rho)}$, where $\epsilon>0$ is chosen
so that $\sigma>2(1+\epsilon)\rho$, and let $E_j=I_K(A_j)\cap
f^{-1}(I_K(\exp A_j^{\phi(A_j)}))$. Then
$${\mathcal A}_K(E_j,1)\leq\exp A_j^{\phi(A_j)}=
\exp n_j^{\phi(A_j)/((1+\epsilon)\rho)},$$ and for $j$
sufficiently large
$$Z_{n_j}(A_j)\leq an_j^2\log3A_j=aA_j^{2(1+\epsilon)\rho}\log3A_j.$$
Since $r_{n_j}=n_j^{1/\phi(r_{n_j})}$ and
$\phi(r_{n_j})\rightarrow\rho$, we have $A_j<r_{n_j}/(8e^4)$ for
all $j$ sufficiently large, hence
$$\frac{n_j}{4}\log\frac{r_{n_j}}{4e^4A_j}>\frac{n_j}{4}\log2.$$
\par As $\phi(A_j)/((1+\epsilon)\rho)\rightarrow1/(1+\epsilon)$ as $j\rightarrow\infty$, we conclude
that for all $j$ sufficiently large we have
$$C_K{\mathcal A}_K^{2\sigma}(E_j,1)<\exp\left(\frac{n_j}{4}\log\frac{r_{n_j}}{4e^4A_j}\right).$$
Proposition \ref{P:geap} implies that $|E_j|\leq Z_{n_j}(A_j)$, so
$$A_j^{-\sigma}|E_j|\leq aA_j^{2(1+\epsilon)\rho-\sigma}\log3A_j\rightarrow0$$
as $j\rightarrow\infty$. The theorem follows by Lemma \ref{L:np},
as $|I_K(A_j)|\geq c'A_j^\sigma$, with a constant $c'$ depending
only on $K$.\end{pf}
\par The theorems of Polya and Gelfond state that if an entire
transcendental function takes integer values at all integer
points, or Gaussian integer values at all Gaussian integer points,
then its order is at least 1, and respectively at least 2. Using
Theorem \ref{T:proportion}, we can obtain asymptotic estimates for
the number of integer (or Gaussian integer) points in the disk of
radius $A$, where a function $f$ takes integer (respectively
Gaussian integer) values.
\begin{Corollary}\label{C:qf} Let $K$ be either
${\Bbb Q}$ or ${\Bbb Q}(i\sqrt{p})$, where $p>0$ is a square free
integer. If $f$ is an entire function of order $0<\rho<\sigma/2$
then
$$\liminf_{A\rightarrow\infty}
\frac{|I_K\cap f^{-1}(I_K)\cap\Delta_A|}{|I_K\cap\Delta_A|}=0.$$
\end{Corollary}
\begin{pf} Note that for $z\in K$ we have $\|z\|=|z|$, so
$$I_K\cap\Delta_A=I_K(A),\;I_K\cap f^{-1}(I_K)\cap\Delta_A=
I_K(A)\cap f^{-1}(I_K(\exp A^{\phi(A)})),$$ for every $A>0$. The
conclusion now follows from Theorem \ref{T:proportion}.\end{pf}
\par We conclude by considering entire transcendental
functions which have a covering system of admissible intervals
$I(R_j,\alpha,\beta,\gamma, C)$ (see Corollary \ref{C:me2}).
Classes of such functions were constructed in Section
\ref{S:scof}. In this case we can estimate ${\mathbf a}_K(s,r,m)$
for fixed values of $r,m$ and for all $s$ sufficiently large. Let
$\tau=1+1/\gamma$.
\begin{Theorem}\label{T:algineq} Let $f$ be as above and let $m,r\geq1$. There
exist positive constants $a$ depending only on $f$, and $C'_K$
depending only on $K$, such that
$${\mathbf a}_K(s,r,m)\ge
\frac {(ms)^{1/\tau}}{64\sigma\rho\tau
a^{1/\tau}}\,\log\frac{ms}{a}-C'_K,$$ for all $s$ sufficiently
large.\end{Theorem}
\begin{pf} By Corollary \ref{C:nz} there is $n_r$ such that
$$Z_n(r)\le an^\tau,\;a=10C(2\beta^{-1})^{1/\gamma},$$ when $n\ge n_r$.
We fix $n_0=n_0(m,r)\geq n_r$ such that
$$\left(m(n+1)^{2\sigma-1}\right)^{-m/n}\ge1/2,\;4e^4r\leq n^{1/(4\rho)},\;
r_n\geq n^{1/(2\rho)},$$ for $n\ge n_0$.
\par Let $s>a(2n_0)^\tau/m$, and let $E$ be a subset of $\Delta_r$ with
$|E|\geq s$. If
$$n=\left[\left(\frac{ms}a\right)^{1/\tau}\right]-1,$$
then $n>n_0$ and $m|E|>Z_n(r)$. Applying Theorem \ref{T:geap} as
in the proof of Theorem \ref {T:id}, it follows that
$$2C_K{\mathcal A}_K^{2\sigma}(E,m)\geq\exp\left(\frac n4\log\frac{r_n}{4e^4r}\right)
\geq\exp\left(\frac{n\log n}{16\rho}\right).$$ Since
$$n\log n\ge\frac{1}{2\tau}\left(\frac{ms}{a}\right)^{1/\tau}
\log\frac{ms}{a}\,,$$ we obtain
$$2C_K{\mathcal A}_K^{2\sigma}(E,m)\ge
\exp\left(\frac{1}{32\rho\tau}\left(\frac{ms}{a}\right)^{1/\tau}
\log\frac{ms}{a}\right),$$ so $${\mathbf
a}_K(s,r,m)\ge\frac{(ms)^{1/\tau}}{64\sigma\rho\tau a^{1/\tau}}\,
\log\frac{ms}{a}-\frac{\log(2C_K)}{2\sigma}\,.$$
\end{pf}
\par We remark that versions of Theorem \ref{T:proportion} and
Corollary \ref{C:qf} can be stated for functions $f$ as in Theorem
\ref{T:algineq}, by requiring that $\rho<\sigma/\tau$ and
replacing the $``\liminf"$ in the conclusion by $``\lim"$.


\begin{thebibliography}{XXXX}
\bibitem[B]{B} A. Baker, {\em Transcendental Number Theory},
Cambridge University Press, 1975
\bibitem[BBLT]{BBLT} L. Bos, A. Brudnyi, N. Levenberg, V. Totik,
{\em Tangential Markov inequalities on transcendental curves},
Constr. Approx., {\bf 19} (2003), 339--354
\bibitem[Br]{Br} A. Brudnyi, {\em Local inequalities for
plurisubharmonic functions,} Ann. of Math. (2), {\bf 149} (1999),
511--533
\bibitem[CP1]{CP1} D. Coman, E. A. Poletsky, {\em Bernstein-Walsh
inequalities and the exponential curve in ${\Bbb C}^2$}, Proc.
Amer. Math. Soc., {\bf 131} (2003), 879--887
\bibitem[CP2]{CP2} D. Coman, E. A. Poletsky, {\em Measures of
transcendency for entire functions}, Mich. Math. J., {\bf 51}
(2003), 575-591
\bibitem[D]{D} J. Dufresnoy, {\em Sur les domaines couverts par les
valeurs d'une fonction m\'eromorphe ou alg\'ebroide}, Ann. Sci.
\'Ecole Norm. Sup., {\bf 58} (1941), 179--259
\bibitem[FN]{FN}Ch. Fefferman, R. Narasimhan, {\em On the
polynomial-like behavior of certain algebraic functions,} Ann.
Inst. Fourier (Grenoble), {\bf 44} (1994),  1091--1179
\bibitem[G]{G} A. O. Gelfond, {\em Algebraic and Transcendental
Numbers}, Dover, 1960
\bibitem[GS]{GS} F. Gramain, F. J. Schnitzer, {\em Ganze
ganzwertige Funktionen: Historische Bemerkungen,}  Complex methods
on partial differential equations,  Math. Res., {\bf 53},
151--177, Akademie-Verlag, Berlin, 1989
\bibitem[H]{H} W. K. Hayman, {\em Meromorphic functions}, Oxford,
1964
\bibitem[K]{K} M. Klimek, {\em Pluripotential theory},
Oxford Univ. Press, New York, 1991
\bibitem[La]{La} S. Lang, {\em Introduction to Transcendental
Numbers,} Addison-Wesley, 1966
\bibitem[Le]{Le} B. Ya. Levin, {\em Distribution of Zeros of Entire
Functions}, Transl. Math. Monographs, vol. 5, Amer. Math. Soc.,
Providence, RI, 1964
\bibitem[M]{M} K. Mahler, {\em Lectures on Transcendental
Numbers}, Lect. Notes Math., {\bf 546}, Springer, 1976
\bibitem[RY]{RY} N. Roytwarf, Y. Yomdin, {\em Bernstein classes},
Ann. Inst. Fourier (Grenoble), {\bf 47} (1997), 825--858
\bibitem[Sa]{Sa} A. Sadullaev, {\em An estimate for polynomials on
analytic sets}, Math. USSR Izvestiya, {\bf 20} (1983), 493--502
\bibitem[Sc]{Sc} T. Schneider, {\em Ein Satz \"uber ganzwertige
Funktionen als Prinzip f\"ur Transzendenzbeweise,} Math. Ann.,
{\bf 121} (1949), 131--140
\bibitem[So]{So} M. L. Sodin, {\em Zeros and units of entire
functions}, Ukr. Math. J., {\bf 40} (1988), 91--95
\bibitem[St]{St} E. G. Strauss, {\em On entire functions with
algebraic derivatives at certain algebraic points,} Ann. Math.,
{\bf 52}, (1950), 188--198
\bibitem[Ti1]{Ti1} R. Tijdeman, {\em On the number of zeros of general
exponential polynomials}, Indag. Math., {\bf 33} (1971), 1--7
\bibitem[Ti2]{Ti2} R. Tijdeman, {\em On the algebraic independence of
certain numbers}, Indag. Math., {\bf 33} (1971), 146--162
\bibitem[T]{T} E. C. Titchmarsh, {\em The Theory of the Riemann
Zeta-Function}, Oxford, 1986
\bibitem[W]{W} M. Waldschmidt, {\em P\'olya's theorem by Schneider's
method,} Acta Math. Acad. Sci. Hungar., {\bf 31}  (1978), 21--25
\end{thebibliography}
\end{document}